%Final variant Apr.30, 2004.

%\magnification=\magstep1
%\baselineskip=12pt
%\parskip=6pt
%\documentstyle{amsppt}
%\voffset=-3pc
%\loadmsbm
\input amstex
\magnification=\magstephalf
\documentstyle{amsppt}
\parindent 20 pt
\NoBlackBoxes
\def\bH{\Bbb H}
\def\ds{\displaystyle}
\def\bR{\Bbb R}
\def\bC{\Bbb C}
\def\bZ{\Bbb Z}
\def\bT{\Bbb T}

\def\sP{\Cal P}
\def\sF{\Cal F}

\def\sH{\Cal H}
\def\sD{\Cal D}
\def\r{\rho}

\def\la{\langle}
\def\ra{\rangle}
\def\df{\dsize\frac}
\def\var{\varepsilon}
\def\Spec{\text{Spec }}

\def\z{\zeta}

\define \hK {$\hat K$}

\define \a{\alpha}

\define \dl{\delta}
\define \g{\gamma}

\define \lm{\lambda}

\define \om{\omega}
\define \s{\sigma}

\define \ve{\varepsilon}

\define \iy{\infty}
\define \p{\partial}

\define \ri{\rightarrow}

\define \sbt{\subset}

\define \edm{\enddemo}
\define \ep{\endproclaim}

\define \1{^{-1}}
\define \2{^{-2}}

\define \BC{\Bbb C}

\define \BR{\Bbb R}
\define \BT{\Bbb T}
\define \BZ{\Bbb Z}

\define \liml{\lim\limits}
\define \intl{\int\limits}

\def\eps{\epsilon}

\define \cp{\operatorname{Cap}}

\nologo
\baselineskip 20pt
\topmatter
\title{A generalization of trigonometric convexity and
its relation to
 positive harmonic functions in
homogeneous domains}\endtitle
\rightheadtext{Trigonometric
convexity}

\author V.~Azarin\qquad D.~Drasin\qquad
P.~Poggi--Corradini\footnote"*"{Supported in part by the National
Science Foundation, Grant No. 
9896337 and No. 9706408\hfill\break}\endauthor
\abstract{We consider functions which are subfunctions with
respect to the differential operator
$$L_\rho = \frac{\p^2}{\p x^2} + \frac{\p^2}{\p y^2} + 2\rho
\frac{\p}{\p x} + \rho^2 $$ and are doubly periodic in $\BC.$ \quad
These functions play an important role in describing the asymptotic
behavior of entire and subharmonic functions of finite order \cite
{6, Ch.3}.  \quad In studying their properties we are led to problems
concerning the uniqueness of Martin functions and the critical
value for the parameter $\r$ in the homogeneous boundary problem for the
operator $L_\r$ in a domain on the torus.}  \endabstract \endtopmatter
\baselineskip 20pt \vskip .10in

\document \subheading{Introduction} The relation between entire
functions and potential theory is a long-standing theme in complex
analysis.  

Consider an entire function $f$ of order $\rho,\,
0<\rho<\infty,$ mean type, i.e., letting
$M(r,f)=\max_\theta\log|f(re^{i\theta})|,$ we have that
$0<\limsup_{r\to\infty}r^{-\rho} M(r)<\infty.$ The classical
Phragm\'en-Lindel\"of indicator corresponding to $f$ is defined as
$$h(\theta) (=h_f(\theta))
=\limsup_{r\to\infty}\frac{\log|f(re^{i\theta})|}{r^\rho},$$
and its key property (cf.\cite {16}) is  that it
is
$2\pi$-periodic and (in the sense of
distributions)
$$h''(\theta)+\rho^2h(\theta)=\nu(d\theta)\geq 0\qquad(0\leq \theta \leq 2\pi),\tag0.1$$
where $\nu$ is  a (positive) measure.  A $2\pi$-periodic  function $h$
 which satisfies (0.1)
is called
{\sl $\rho $-trigonometrically convex} ($\r$--t.c.).
The behavior of solutions to (0.1)  reveals many facts about entire
functions of finite order
(\cite {16}).
Note that  (0.1)
holds if and only if
$$v(z)=r^\rho h(\theta)\tag0.2$$
is subharmonic in the plane.  To see the connection between
(0.2) and entire functions of order $\rho$, mean type, recall that $f$ is
 of completely regular
growth if for some $0<\rho <\infty\ $ the following
limit exists
$$D'-\lim_{t\ri\infty}\log |f(zt)|
t^{-\rho},\tag 0.3$$
in the distributional topology $D'(\Bbb C\setminus 0)$.
Calling this limit $v (z)$ it will then  necessarily
have the form (0.2) (see \cite {16},\cite{6, Ch.3}).

This class of functions
was introduced (in a modified form) independently by B. Ja. Levin
and A. Pfluger, and is
a major focus of \cite {16}.
A more complicated asymptotic behavior arises from the class of entire
functions with {\sl periodic limit set} \cite {6, Ch. 3}.  This class
is based on a $T$-{\sl automor\-phic} subharmonic function $v(z)$;
i.e., a subharmonic function for which there are fixed $T>1$, $\r>0$
such that $$v(Tz) = T^\r\,v(z)\quad(z\in \Bbb C).\tag 0.4$$
Given such a function $v$, we thus consider
an entire function $f$ which satisfies a condition analogous to (0.3):
for every $1<\tau \leq T$ there exists a sequence $t_j\ri\infty$ such
that $$D'-\lim_{t_j\ri\infty}\log |f(zt_j)| t_j^{-\rho}=v(z\tau)\tau
^{-\rho}$$
 and for every sequence $u_{t_j}(z)=\log |f(zt_j)| t_j^{-\rho}$ there 
exists a
$\tau $ as above and a subsequence converging to $v (z\tau)\tau^{-\rho}$.
If $v$ happens to satisfy (0.4) for every
$T>0$, as in (0.2), then
we recover (0.3), because the family
$\{u_t(z):
t\in [1,\infty)\}$
is always compact in $D'(\Bbb C\setminus 0).$

The functional equation (0.4) makes it natural to consider functions
defined on open sets $G$ which are invariant under multiplication by
$T,$ i.e., $TG=G$. We call such $G$ a {\sl $T$-homogeneous set}, and
reserve the notion {\sl $T$-homogeneous domain} to indicate that $G$
is open and connected. The boundary of a nonempty $T$-homogeneous set $G$
(not $\bC$)
always includes $0$ and $\infty,$ and 
we always assume that $\partial G$ (and $\partial D$, below)
has positive capacity.

Consider the class $\sP$ of positive harmonic functions on
$T$-homo\-gen\-eous domains $G$ which are bounded in any bounded
subset of $G$ and which vanish quasi-everywhere (i.e., outside a set
of zero capacity) on $\partial G$ \cite {3, 10, 14, 18}. For a general
$T$-homogeneous domain, the class $\Cal P$ may contain infinitely many
non-proportional functions
(see an example in \S 3.)
We identify a subclass $\Cal F\subset
\Cal P$ consisting of functions of restricted growth at infinity, as
in (0.6) below, which turns out to be always non-empty and
one-dimensional: it consists of positive multiples of a single
function.  We show in \S 5 that $\Cal P = \Cal F$ for a large class of
domains.  \proclaim {Theorem 0.5} Let $G\subset \Bbb C$ be a
$T$--homogeneous domain.  Let the family $\sF\subset \Cal P$ consist
of functions $v\in \Cal P$ such that
$$
M(r,v)\equiv \max\limits_{|z|=r\atop z\in G} v(z)\leq C r^k\quad(r>r_0)
\tag 0.6$$
for some $r_0=r_0(v)<\infty$ and $0\leq k=k(v) < \infty$ $($i. e., $v$
has finite order$)$. 
Choose some $z_0 \in G,\ |z_0|=1$.  Then
\roster
\item there exists a unique function
$H \in \Cal F$ with $H(z_0)=1$ and hence
$v\in\sF\Leftrightarrow v=cH$ for some constant $c>0$;
\item there exists a unique
$\r (G)>0$ such that every $v\in \Cal F$
satisfies the functional equation
$$v(Tz)=T^{\r(G)} v(z) \tag 0.7$$
\endroster

\endproclaim

Let us note that the equation (0.7)
coincides with (0.4) for
$\r=\r(G)$.

Many properties of $G$ are reflected in $\rho(G),$ and will be discussed
in, for example, Theorem 0.17, \S 4, and \S 6.
In \S 4.6 we present several interpretations of 
$\r(G)$ when $G$ is simply-connected and $T$-invariant.

Now let $v$ satisfy (0.4). Then the function
$$  q(z) = v(e^z)e^{-\rho x}\tag0.8$$
is $2\pi$--periodic in $y$ and periodic in $x$ with
period $P=\log\, T.$
The function $q$ can be considered as a function on  a torus $\BT^2_P,$
 obtained by identifying the opposite sides of the rectangle
$R = (0,P) \times (-\pi, \pi).$
The homology
group of $\BT^2_P $ is nontrivial, with basis the cycles
$\gamma_1,\gamma'_1$, where $\gamma_1=\BT^2_P\cap \{y=0\},\
\gamma'_1=\BT^2_P\cap \{x=0\}$.

Let $\pi$ be the covering map of $\Bbb C$ onto $\Bbb T^2_P$, then
$\phi =\pi \circ \log$ is a well-defined covering map of $\bC
\setminus\{0 \}$ onto 
$\Bbb T^{2}_{P}$, where the group of deck transformations is given by
the dilations by $T^{m}$ for $m\in \bZ$. 
So if $G$ is a given $T$-homogeneous domain, then
$$D = \pi \circ \log G=\phi (G),\tag 0.9$$
is a domain in $\Bbb T^2_P$.
On the other hand, not every domain in $\BT^2_P$ has a
$T$-homogeneous domain as its preimage under $\phi$. 
The preimage
$\phi^{-1} (D)$ under $\phi$ is a possibly disconnected set which is invariant
under dilations by $T^{m}$ for $m\in \bZ$.
An intrinsic description is given by the
next proposition.

 \proclaim {Proposition 0.10} Let $\hat\g$ be a closed curve in a
domain $D\sbt\BT^2_P$ homologous in $\BT^2_P$ to a cycle $\gamma =
n_1\gamma_1+n'_1\gamma'_1,\ n_1,n_1'\in \BZ.$ Then

1. If no such $\hat\g$ can be found in 
$D$ so that $n_1\neq 0$, then 
$$\phi^{-1}(D)=\cup_{j=-\iy}^{\iy}G_j,$$ where $G_j=T^jG_0,$ $G_0$ is
an arbitrary connected component of $\phi^{-1}(D),$ and $G_j\cap
G_l=\emptyset$ for $j\neq l.$
 
2. If there exists a curve $\hat\g$ as 
above with $n_1\neq 0,$ then
$$\phi^{-1}(D)=\cup_{q=0}^{k-1}G_q,$$ where $k =\min |n_1|$ with the
minimum taken over all such curves $\hat\g$; $G_0$ is an arbitrary
component of 
$\phi^{-1}(D);$ $G_j, \ j=0,1,...,k-1$,
are disjoint $T^k$-homogeneous domains, and for every $m\in\BZ,\ T^{
m} G_0=G_q$, 
provided $ \ m=lk+q,$ for some 
$q\in \BZ,\ 0\leq q\leq k-1,\ l\in \BZ.$ \ep

We call domains as in part 2 of Proposition 
0.10 {\it connected on spirals}. In particular, this proposition 
shows that  for every $D$ connected on spirals, we can  find a connected $T^k$-
homogeneous domain that relates to $D$ by (0.9).
The proof of Proposition 0.10 is given in \S 2.

Let us give examples. The domain $D^\prime=\BT^2_P\cap \{|x-P/2|<P/4\}$ is not
connected on spirals, while $D^{\prime\prime}=\BT^2_P\cap \{|y| <
\pi/4\}$ is.  It follows that $D^\prime\cap D^{\prime\prime}$ is not
connected on spirals while $D^\prime\cup D^{\prime\prime}$ is.

The situation can be more complicated.
Let $R$ be the rectangle
$[0,P]\times[-\pi, \pi]$.  We construct a network of disjoint strips
$\{D_j\}_{-\infty}^{+\infty}$ which connect the vertical portions of
$\p R$.  If
$I_j = D_j\cap \{x = 0\}$ and $I'_j = D_j\cap\{x = P\}$,
we arrange that $I_j$ and $I'_{j-1}$ have the same projection
on the $y$-axis, and these projections cluster to
$\{y= \pm \pi \}$ as $j \to \pm \infty$.
After the usual identification of sides of $\p R$, we obtain
a domain $D\subset \BT^2_P$ which winds infinitely often in the
 $x$-direction on $\BT^2_P$.  This $D$ is not connected on spirals
and illustrates how the "bad" case in the proof of Proposition 2.1 looks
like. More complicated domains yet are obtained by replacing the
fundamental rectangle $R$ by a fundamental parallelogram $R'$
whose vertical sides are  $\{x=0,-\pi<y<\pi\}$ and $\{x=P, -\pi+k2\pi< y
<\pi+k2\pi\}$ for some fixed integer $k$. Then, the strips $D'_j$ are in
$R'$ and connect $I_j$ to $I''_j$, but now the projections of $I_j$ and
$I''_{j-1}$ on the $y$-axis differ by a translation of $k2\pi$ units.

One more example. Consider the family of lines
$L_l:=\{z=x+iy:y=\pi/(kP)x+l\pi/k,\ x\in \BR\}, \ l\in \BZ.$ It
determines a closed curve (spiral) $\hat\g$ on $\BT^2_P$ with $n_1=k.$
The open set $D_k=\{z:|z-\z|<\eps, \z\in L_l,\ l\in\BZ\},\
0<\eps<P/2\sqrt {\pi^2+k^2},$ determines a domain $\hat D_k$ on
$\BT^2_P$ that is connected on spirals, and such that $\phi^{-1}(\hat
D_k)$ consists of $k$ components, each $T^k$-homogeneous.  \subheading
{0.11} Since the function $v$ of (0.4) is subharmonic, the function
$q$ of (0.8) is upper semicontinuous and in the $D'$ topology on $\Bbb
T^2_P$ satisfies the inequality $L_\r q\geq 0$, where
$$
L_\rho:=\Delta+2\rho {\partial\over \partial x}+\rho^2.\tag0.12
$$
That is, $L_\r q $ is a positive
measure on $\BT^2_P.$

The operator $L_\r$ arises naturally since
(0.8) shows that if $v (Z)$ is a smooth function and $q$
is related to $v$ as in (0.8), then
$$\Delta_{Z}v(Z)=e^{(\rho-2) x}L_\rho q(z),\quad
 Z=e^z,\ z=x+iy.\tag0.13$$
Such functions $q$ are called {\sl subfunctions with respect to }
$L_\r,$ or $L_\r$-{\sl subfunctions}.  Note that $L_\rho $ is not
symmetric when $\rho \neq 0$.

In this paper we obtain some properties of $L_{\rho}$-subfunctions;
these generalize those of $\r$-t.c.functions.  For the
theory of entire functions modeled on functions $v(z)$ as in (0.4), the
$L_{\rho}$-subfunctions play the same role that the $\r$--t.c.functions play
for entire functions of completely regular growth (see \cite {2, 3,
4}).  \subheading{0.14} The study of $L_{\rho}$-subfunctions depends
on properties of the operator $L_{\r}$ for arbitrary $\r$ and, in
particular, on properties of solutions to the homogeneous boundary
problem
$$\aligned
&L_\r q = 0\,  \quad \text{in }D;\\
&q\bigm|_{\p D} = 0,\endaligned \tag0.15$$
where $D$ is a domain in ${\BT}^2_P$
and $q$ is bounded in $\p D$ with
boundary value zero quasi-everywhere.
This is a spectral problem for a {\sl pencil
} of differential operators (the standard reference is \cite {17};
cf. \S 1 below).

We emphasize that in principle a solution of this problem can be
defined for an arbitrary domain $D\subset{\BT}^2_P$; recall, however,
that the boundaries of all domains considered here have positive
capacity.
  
The {\sl spectrum} of the problem (0.15) consists of those (complex) $\rho$ for
which (0.15) holds for some function $q\not\equiv 0$.  
We identify when the spectrum is nonempty,  and give some basic properties in
Propositions  1.36, 1.37.  We also show that the {\sl minimum positive point}
of this spectrum, $\rho(D)$, exists, and is intimately connected with
the function $H$ and the number $\rho (G)$ produced in Theorem 0.5,
with some component  $G$ of $\phi^{-1}(D)$ 
(see (0.9)).

\proclaim {Theorem 0.16}The
following hold: \roster
\item $\r(D)<\iy$ iff $D$ is connected on spirals;
\item If $\r(D)<\infty$, then $\r(D)=\r(G)$,
and up to a constant multiple the corresponding
eigenfunction is
$$q(z)=H(e^z)e^{-\r(D)x}.$$
\endroster
\endproclaim

A property of $\r (D)$ which carries over from the classical potential
theory is {\sl strict monotonicity}: \proclaim {Theorem 0.17} Let
$D_1,\ D_2$ be domains on $\Bbb T^2_P$ which are connected on spirals.
If $D_1\sbt D_2 $ and $\cp\ (D_2\setminus D_1)>0$ then the strict
inequality $\r(D_1)>\r(D_2)$ holds.  \endproclaim

\subheading{0.18}
In \S 1 we find the fundamental solution (for $\r\notin \BZ$) and 
the generalized fundamental solution (for $\r \in \BZ$) 
of the equation 
$L_\r q(z)=0$ on the whole torus $\BT^2_P.$ In \S 8.14
we use it in a representation    
 which is a generalization of the well known representation of
$\r$--t.c. functions from \cite {16, Theorem 24} (see also  
\cite {6, \S 2,(4),(5)}). For application of 
this representation, see \cite {2}.

 In \S 7 we introduce the Green function for
$L_\rho$, and use this as basis for studying $L_\rho$-subfunctions on
subdomains of the torus $\BT^2_P$.

We also consider subharmonic minorants of a given real function $m$ in
the plane (see, for example, \cite{14}). In application to subharmonic
functions with periodic limit sets this leads us to considering of
$L_\r - subminorants\ of\ a\ function\ m,$ i. e., $L_\r$--subfunctions
$u(z)$ with $u(z) \leq m(z)$ for $z\in\BT^2_P.$ Theorems 9.15 and 9.16
imply

\proclaim{Theorem 0.19} Let $m$ be a continuous function on
$\BT^2_P$.  If $m$ has a non-zero $L_\rho$ subminorant, then  
$\rho(D)\leq\rho$ for some component $D$ of the set
$\Cal M_+:=\{z: m(z)>0\}$. 
 
Conversely, if $\rho(D)<\rho$ (strict inequality!) for some component
$D$ of the set $\Cal M_+$, and $m(z)\geq 0$ for all $z\in\BT^2_P,$
then $m$ has a non-zero $L_\rho$-subminorant.\endproclaim

This generalizes properties of $\rho$-t.c. functions,
since when  $D=\{z\in \BT^2_P: \Im z \in (\alpha, \beta)\}$
we have $\rho(D)=\pi/(\beta - \alpha).$

The borderline case $\rho(D)=\rho$ depends essentially on the behavior of
$m$ near $\partial D$, and warrants further scrutiny,  
as well as the case when $m$ changes its sign.

There is a specific question that arises in studying the completeness of
exponential systems \cite{4}.
An $L_\rho $--subfunction  $u(z),\ z\in \BT^2_P
$ is {\it minimal} if the
function $m(z) = u(z)-\epsilon $ does not have an $L_\rho $--subminorant in
$\BT^2_P$ for arbitrarily small $\epsilon > 0.$
A full description of minimal functions is not known (see \cite {8,
Problem 16.9}), but some necessary and some sufficient
conditions are obtained here.   For example, in \S 9 we show
\proclaim{Theorem 0.20} Let $\Cal H_\rho(u)$ be the maximal open set on which
$L_\rho u=0$.  If there exists a connected component
$M\subset \Cal H_\rho(u)$ such that
$\rho(M)<\rho,$ then $u$ is a minimal
$L_\rho$-subfunction. \endproclaim

This paper is organized as follows:

In \S 1 we study properties of the operator $L_\r$ and the generalized boundary
problem,
and prove  Theorem 0.5  in \S 3.
Theorem 0.16 and 0.17 are proved in \S 6, and \S\S 7-9 are devoted to
$L_\rho$-subfunctions and subminorants.

We are grateful to Profs. A. Ancona, A. Marcus,  M. Sodin and especially to
Profs.  A. Eremenko and V. Matsaev for very valuable discussions and
suggestions.

\subhead 1.\ The operator $L_\rho$; characterization of $\Spec
D$\endsubhead First we study fundamental solutions of $L_\rho$ on 
the whole torus
 $\BT^2_P$.  \proclaim{Proposition 1.1} If $\rho \notin \Bbb Z=\{0,
\pm 1, \pm 2,\dots\}$, the operator $L_\rho$ has a unique fundamental
solution $E_\rho(z)$ on $\BT^2_P$ with singularity at $0\ (+kP+2\pi
li,\ k,l\in \BZ) $.  \endproclaim \demo{Proof} We solve the equation
$L_\rho E_\rho = \delta_0(z)$ in $\Cal D'(\Bbb T^2_P)$ (the space of
distributions), where $\delta_0$ is the Dirac function supported at $0
\in \Bbb T^2_P$.  Using the transformation $x'= 2\pi x/P,\ y'=y,$ we
obtain that the period in $x'$ is $2\pi$, so that
$\bT^2_P\equiv\bT^2_{2\pi}\equiv \bT^2,$ and $L_\rho$ is replaced by
$$L_{\rho,P} = \displaystyle \left(\frac{2\pi}{P}\right)^2\frac{\p^2\
}{\p x^2}+\frac{\p^2\ }{\p y^2} +2\rho\left(\frac{2\pi}{P}\right)
\frac{\p\ }{\p x} +\rho^2.\tag1.2$$ The Dirac function is
characterized by the action $\la \delta_0, g\ra=g(0)$ for $\
g\in \sD(\bT^2)$, the class of infinitely-differentiable doubly
$2\pi$-periodic functions.  The system
$\phi_{k\ell}(z)=e^{ikx}e^{i\ell y}\,(k,\ell \in \bZ)$ is dense in
$\sD(\bT^2)$.  We compute the Fourier coefficients of the solution
$E,$ corresponding to (1.2) by solving
$$
\la L_{\rho}E,\phi_{k\ell}\ra \equiv \la
E, L^*_{\rho}\phi_{k,\ell}\ra=1,\qquad(k, \ell \in \bZ),\tag1.3
$$
where $L^*_{\rho,P} = L_{-\rho,P}$ is 
symmetric to 
 $L_{\rho,P}$.   Equations (1.2) and (1.3)  then yield that
$$
a_{k\ell}\equiv \la E, \phi_{k\ell}\ra
= \left[-\left(\frac{2\pi}{P}\right)^2 k^2-\ell^2+2\rho\frac{2\pi}{P}ik
+\rho^2\right]^{-1},\tag1.4
$$
and since $\rho \notin \bZ$, these coefficients are uniquely determined.

Note that if $E_1$ and $E_2$ are solutions, then all Fourier
coefficients of $E_1-E_2$ vanish.  Thus $E$ is unique, and $$ E(z) =
\sum_{k,\ell} a_{k\ell}e^{ikx}\cdot e^{i\ell y}, $$ where $a_{k\ell}$
are determined by (1.4).  The series defining $E$ always converges in
the sense of distributions, and it is well-known that solutions to
elliptic homogeneous equations with constant coefficients are
real-analytic (cf.  \cite {11}). Thus $E$ is smooth when $z\neq 0$ and
has logarithmic singularity at $z=0.$
 
Finally, we set
$$E_\r (z):=E\left(\frac {xP}{2\pi}+iy\right).$$
\qed \enddemo

For nonintegral $\r$ we may express the fundamental solution in another form,
using common notions from the theory of subharmonic
functions, which also
provides an independent  way to check the regularity of
$E_\rho$ off the diagonal.
Let $p = [\r]$ and $H(u,p)$ be the logarithm of the classical Weierstrass factor of genus $p:$
$$
H(u,p) = \log |1-u|+ \Re\left(\sum^p_{k=1} \df{u^k}{k}\right),\tag1.5
$$
and set
$$
g(z,p,\r) = H(e^z, p)\ e^{-\r x};\qquad
g(z,p,\r,P) = \sum^\infty_{k=-\infty} g(z+kP, p,\r).\tag 1.6
$$
This series converges for all $z\neq kP,\ k\in \bZ$ because of the
inequalities (see, for example, \cite {9, Lemma 1.5})
$$
|H(u,p)|\leq\cases
 C_{p,\varepsilon}|u|^p,& \text{for}\ |u|\geq 1+\varepsilon\\
C_{p,\varepsilon}|u|^{p+1},&\text{for}\ |u|\leq 1-\varepsilon,\\
 C_{p,\varepsilon},
&\text{for}\  1-\varepsilon\leq |u|\leq 1+\varepsilon, |u-1|>\varepsilon,\endcases
\tag1.7$$
where $C_{p,\varepsilon}$  is a constant independent of $u$.
Thus $g(z,p,\r,P)$ is $2\pi$-periodic in $y$ and
$P$-pe\-riodic in $x,$
and so may be viewed as a function on $\BT^2_P$.

\proclaim{Proposition 1.8}Let
$g(z,p,\r,P)$
 be from $(1.6)$ and $E_\r$ be the
fundamental solution given by Proposition 1.1.
Then $$
\df{1}{2\pi}\ g\left(z, p, \r, P\right) = E_\r(z).  $$
\endproclaim
\demo{Proof} We show that $({1}/{2\pi})g$ is a fundamental solution of $L_{\r}$.
To prove this we can use test  functions from
$\Cal D(\BT^2_P)$ supported only near one
singularity of  the
series in (1.6). We can suppose also that it is
the term for $k=0.$

Since all other terms satisfy (0.15) , we need only show that
$$L_\r g(\cdot,p,\r)=2\pi\dl_0 \tag 1.9$$
on functions
$F\in \Cal D((-P/2,P/2)\times(-\pi,\pi)).$  Our arguments use the
correspondence (0.8).

To avoid confusion, we take  $z$ for the variable on
$\BT^2_P$, and $Z$ for regions in the plane (here $Z=e^z$).
 Thus let $F$ have compact support near $z=0$ viewed
as a point in $\BT^2_P$.  Then using the calculus of distributions with (0.13) and (1.6)
we have
$$
\aligned
I:=&<L_\r g(\cdot,p,\r),F>=<g(\cdot,p,\r),L_{-\r}F>=\int L_{-\r}F(z)g(z,p,\r)dxdy=\\
&=\int[e^{(-\r-2)x} L_{-\r}F(z)][e^{\r x}g(z,p,\r)]e^{2x}dxdy\\=&\int e^{(-\r-2)x}
L_{-\r}F(z)H(e^z)e^{2x}dxdy
=\int \Delta_Z(F(z)e^{-\r x}) H(Z)dXdY.
\endaligned$$
However, $(1/2\pi)\log|z-1|$ is the fundamental solution to the  Laplace equation
at $z=1$, and so if we set $f(Z)= F(z)e^{-\r x}$ (recall  (0.8)),
we find that $I = 2\pi f(1)=2\pi F(0), $ and this  gives (1.9).

By uniqueness of the fundamental
function on $\BT^2_P$ with singularity at zero,  we obtain  Proposition 1.8.
\qed
\enddemo

When $\r$ is an integer, the  reasoning which gave Proposition 1.1 gives (proof omitted)

\proclaim{Proposition 1.10}\ If $\r\in\bZ$,
there exists a generalized fundamental solution $E'_\r(z)$ satisfying
the equation
$$
L_\r\ E'_\r (z) =  \dl_0 (z) - 2\cos \r y\tag 1.11
$$
on $\bT^2_P$.

The function
$E'_\r (z)$ is defined uniquely up to an addend of the form
$$
A \cos \rho (y-y_0),
$$
where $A$ and $y_0$ are arbitrary.
\endproclaim

We next define the spectrum of the problem (0.15) for a domain with
arbitrary boundary and prove that it is a discrete set without any
finite point of condensation (Proposition 1.36). The natural way for
this is to transform the problem (0.15) to an integral equation.

The domain $D\sbt \BT^2_P$ is a Riemannian space of hyperbolic type
and admits a Green function for the Laplace operator on $\BT^2_P$.
(see, e.g., \cite {1, Ch.10}).

We will use the local coordinates $z=x+iy$ 
which preserves not only the sign of  the Laplace 
operator but also the Laplace operator itself.
Thus these coordinates preserve harmonicity and 
subharmonicity of functions as well as their mass 
distributions. 

We denote the Green function by $g(z,\z,D)$ and 
extend $g$ to 
the whole
$\BT^2_P$ by defining 
$g(z,\z,D)=0$ when $z$ or $\z \in
\bT^2_P\setminus D.$ Denote by $\nabla_z$ the 
gradient operation in $x,y;\ dz=dxdy$ for 
$ z=x+iy$ and $d\z=d\xi d\eta$ for 
$ \z=\xi+i\eta.$
\proclaim{Proposition 1.12} Let $D\subset \bT^2_P$ be a domain.
Then there  exist  constants $C=C(p)$ such that
$$\align &\sup_{z\in \Bbb T^2_P} \int_{\Bbb T^2_P }|g(z,\zeta, D)|^p\,
d\zeta\leq C\qquad (0\leq p <\infty);\\
&\max\left \{\sup_{z\in \Bbb T^2_P}\int_{\Bbb T^2_P}
|\nabla_\z g(z, \zeta, D)|^p d\zeta ,
\sup_{\z\in \Bbb T^2_P}\int_{\Bbb T^2_P}
|\nabla_\z g(z, \zeta, D)|^p dz\right \}
\leq C \quad(0\leq p <2),
  \endalign  $$
\endproclaim

\demo {Proof} The function $-g(z,\z,D)$ is subharmonic in
$\BT^2_P\setminus\{\z\}$ with masses distributed on $\p D$ and the
total mass $1.$ In a neigborhood of $\z$ it is represented in the form
$-g(z,\z,D)=v(z,\z)-\log|z-\z|$ where $v$ is a subharmonic function,
with its masses distributed on $\p D$ and the total mass $1.$ Note, if
the neighborhood does not intersect $\p D,$ then $-g(z,\z,D)$ is
harmonic. As a matter of fact, the mass distribution of $v$ coincides
with the harmonic measure of $D,$ which does exists even for every
Riemannian space of hyperbolic type.

It is possible to check, using the H\"older inequality and the
continuity in $z$ of the functions
$$\int\limits_{U}|\log |z-\z||^pd\z,\ 0\leq p<\iy;\
\int\limits_{U}|z-\z|^pd\z,\ 0\leq p<2,$$ where $U$ is a small disc,
that a potential of bounded masses belongs locally to $L^p,\ 0\leq
p<\iy,$ while its gradient belongs locally to $L^p,\ 0\leq p<2.$ These
are locally true for every subharmonic function, by the Riesz
representation . This also holds on the whole $\BT^2_P$ because of its
compactness.  Thus this is fulfilled for the Green function and its
gradients $\nabla_z$ and $\nabla_\z.$ \qed\edm

Using $\z=\xi + i \eta$ as the local
coordinates on $\bT^2_P$ we  set
$$g^1(z,\z,D):=(\partial g/\partial
\xi)(z,\z,D).\tag 1.13$$  
The functions $g$ and
$g^1$ induce integral operators on
$C^\infty(\Bbb T^2_P)$, which will be the focus of our
attention:
$$
\align
G_Dq(z)&=\int_D g(z,\zeta,D) q(\zeta) d\zeta\equiv \int_{\Bbb T^2_P}
g(z,\zeta,D)q(\zeta)d\zeta,\\
G^1_Dq(z)&=-\int_D g^1(z,\zeta,D) q(\zeta)d\zeta\equiv -\int_{\Bbb T^2_P}
g^1(z,\zeta,D) q(\zeta) d\zeta.
\tag 1.14\endalign
$$
We will use the following properties of these operators: \proclaim
{Proposition 1.15} Let $G_D,\ G^1_D$ be as above.  Then $G_D,G^1_D$
can be extended as compact operators from $L^2(\Bbb T^2_P)$ to
$L^p(\Bbb T^2_P)$ for every $p\geq 1$ (in particular, for $p=2$).  \ep

Before proving this theorem, we present, following \cite {15}, some
preliminary information on integral operator theory, corresponding to
our case.

Set
$$Ku(z):=\int\limits_{\BT^2_P}K(z,\z)u (\z ) d\z.\tag 1.16$$
This is an integral operator with the kernel 
$K(z,\z),\ z,\z\in \BT^2_P.$

Define the functions
$$\phi_r(z):=\left (\int\limits_{\BT^2_P}
|K(z,\z)|^rd\z\right)^{1/r};\tag 1.17$$
$$\psi_{r^*}(\z):=\left (\int\limits_{\BT^2_P}
|K(z,\z)|^{r^*}dz\right)^{1/r^*},\tag 1.18$$
and their norms
$$\|\phi_r\|_q:=\left (\int\limits_{\BT^2_P}
|\phi_r(z)|^qdz\right )^{1/q};\ 
\|\psi_{r^*}\|_{q^*}:=\left (\int\limits_{\BT^2_P}
|\psi_{r^*}(\z)|^{q^*}d\z\right)^{1/q^*}.
\tag 1.19$$
Set also
$$\|u\|_{\a}:=\left (\int\limits_{\BT^2_P}
|u(z)|^{\a}dz\right)^{1/\a};\tag 1.20$$
and define the space $L^{\a}$ as a space 
obtained by the closure of the space of infinitely differentiable functions 
$u(z),\ z\in \BT^2_P$ with respect to this norm. It coincides with the
space of mesurable functions for  
which the integral (1.20) is finite. The space 
$L^{\a}$ is a Banach space for each $\a\geq 1.$
We are going to use the following assertion, 
which is  a restatement of \cite {15, Theorem 7.1}
for our case.
\proclaim {Theorem KZPS}Suppose $K(z,\z)$ satisfies the condition
$$\max\{\|\phi_r\|_q,\|\psi_{r^*}\|_{q^*}\}
<\iy\tag 1.21$$
for some $r,\ q,\ r^*,\ q^*$ such that  
$$\min \{r^*,q^*,r\}\geq 1;q\in (0;\iy).\tag 1.22$$   
Then for every $0<\tau<1,$ the integral operator 
(1.16), acting from 
$L^{1/\alpha(\tau)}$ to $L^{1/\beta(\tau)}$ for
$$ \alpha(\tau):=1 - (1-\tau)(1/r) -\tau(1/ q^*);\ \ 
\beta(\tau):=(1-\tau)(1/q) + \tau(1/ r^*),  \tag 1.23$$
 is compact; and its operator norm satisfies the inequality 
$$\|K\|_{1/\alpha(\tau)\rightarrow1/\beta(\tau)}\leq
\|\phi_r\|_{q}^{1-\tau} \|\psi_{r^*}\|_{q^*}^\tau .\tag 1.24$$
\ep
\proclaim {Proposition 1.25}Set
$$1/r=1/r^*=1/2 + 1/(2p);\ 1/q=1/q^*=1/(2p).
\tag 1.26$$

If (1.21) is satisfied under these conditions, then for every $p\geq
1,$ the operator $K$ acts from $L^{2}$ to $L^{p}$ as a compact
operator; and its norm satisfies the inequality (1.24) for
$\tau=1/p.$ \endproclaim \demo {Proof} The conditions (1.22) are
fulfilled. From (1.23) we have
$$\alpha=1-(1-1/p)(1/2 + 1/(2p))-1/p\  1/(2p)=1/2;
$$
$$\beta =(1-1/p)\ 1/(2p)+ 1/p\ (1/2 + 1/(2p))=1/p.$$
\qed
\enddemo
\demo {Proof of Proposition 1.15} Let us check the 
conditions (1.21). Set in (1.17), (1.18)
$r,q,r^*,q^*$ from (1.26) and 
$K(z,\z):=g(z,\z,D).$ We obtain
$$\phi_r (z)=\left (\intl_{\Bbb T^2_P}|g(z,\z,D)|^rd\z
\right )^{1/r}.$$
Thus $\phi_r (z)\leq C^{1/r},\ z\in \Bbb T^2_P,$ by Proposition 
1.12 (the first inequality). Using (1.19),
we obtain that $\|\phi_r\|_q <\iy$ for every $q.$ Since
$g(z,\z,D)=g(\z,z,D)$ , $r=r^*$ and  
$q=q^*,$ the 
inequality $\|\psi_{r^*}\|_{q^*} <\iy$ also holds.
Thus the compactness of $G_D$ is proved.

From (1.26) $r,r^*<2$ if $p\geq1.$ Set 
$K(z,\z):=-g^1(z,\z,D).$ In (1.17) we have
$$\phi_r(z)=\left (\intl_{\Bbb T^2_P}|g^1(z,\z,D)|^rd\z
\right )^{1/r}\leq \left (\int_{\Bbb T^2_P}
|\nabla_\z g(z, \zeta, D)|^r d\zeta\right )^{1/r}.$$ 
By Proposition 1.12, we obtain 
$\phi_r(z)\leq C^{1/r},\ z\in \Bbb T^2_P.$ Thus 
$\|\phi_r\|_q <\iy$ holds for every $q.$
For proving $\|\psi_{r^*}\|_{q^*} <\iy,$ we set  
$$\psi_{r^*}(\z)=\left (\intl_{\Bbb T^2_P}|g^1(z,\z,D)|^{r^*}dz
\right )^{1/r^*}$$
and repeat the previous reasoning.
\qed
\edm

Now we may rigorously  define
problem (0.15) for an arbitrary domain with no
assumption concerning  its boundary (other than positive capacity).

Consider the {\it operator pencil}
$$\Cal G_D (\r):=I+2\r G^1_{D}+\r^2  G_{D};\tag 1.27$$
i.e., a family of operators acting from
$L^2(\Bbb T^2_P)$
to $L^2(\Bbb T^2_P),\ $ depending  on a complex parameter
$\rho$
(see, e.g., \cite {17}).

Since $q, G_D q, G^1_Dq \in L^2(\Bbb T^2_P)$ they
also belong to $L^1(\Bbb T^2_P)$ and can be
considered as distributions in $D.$

\proclaim {Proposition 1.28}Let $D\sbt
\Bbb T^2_P$ be an arbitrary  domain (whose boundary has positive capacity), and
$q\in L^2(\Bbb T^2_P).$
 If $q$ satisfies the distribution equation
$$\Cal G_D(\r)q=0,\tag 1.29$$
then $q$ is infinitely
differentiable inside $D$,  $L_\r q(z)=0$ for all $z\in D$,
$q$ is bounded, and $q$ tends
to zero at  every regular point of $\p D,$
i.e. $q$ is the solution of (0.15).
\endproclaim

\demo {Proof} Consider the equality (1.29) in $\Cal D'(D)$.  Let
$\phi\in \Cal D(D)$ and apply $\Cal G_D q$ to $({1}/({2\pi}))\Delta
\phi$. The definitions of $G_Dq$ and $G^1_Dq$ and standard properties
of the Laplacian (in particular that $({1}/{2\pi})\Delta
g(\cdot,\z)=\delta_\z)$ imply that $<q,({1}/({2\pi}))\Delta \phi>=<
({1}/({2\pi}))\Delta q, \phi>, <G_D q, ({1}/({2\pi}))\Delta
\phi>=<q,\phi>,$ as well as
$$\align <G^1_D q,\frac {1}{2\pi}\Delta \phi>=&-\int_D<\frac{\p}{\p
\xi}g(\cdot,\z),\frac {1}{2\pi}\Delta \phi>q(\zeta)\,d\z=
-\int_D\frac{\p}{\p \xi}\phi(\zeta) q(\zeta)\,d\z\\ &=<\frac{\p}{\p
\xi}q,\phi>.\endalign$$ Hence $0=<\Cal G_D q,\Delta \phi>= <L_\rho
q,\phi>$. By \cite {11, \S 11} $q$ is infinitely differentiable and
satisfies the same equation in the classical sense, so $L_\r q(z)=0$
holds in $D.$

We next prove that  both integrals which
appear in the operator
$\Cal G_D q$ vanish as  $z\ri z_0  \in D$ for
any regular point
$z_0\in \p D.$

It is clear  that $z_0,$ being regular, ensures that
$\liml_{z\ri z_0}g(z,\z,D)=0.$  Since $g$ is harmonic in $D\setminus
\{\zeta\}$, 
it follows that if $K\sbt D$ is compact, we have uniformly that
$$\liml_{z\ri z_0}g^1(z,\z,D)=0 \ \ (\z\in K).\tag 1.30$$
(recall that $g^{1}$ means derivative with respect to $\zeta$).
We show, for example, that
$$\liml_{z\ri z_0}G^1_Dq(z)=0. \tag 1.31$$
Proposition 1.15 and  the equation
$q=-2\r G^1q-\r ^2 Gq$
show  that $q\in L^p(\Bbb T^2_P)$ for
some $p>2$. Let $\eps$ be arbitrarily small. Choose $K$
such that
$$\left(\intl_{D\setminus K}|q(\z)|^pd\z\right)^{1/p}<(1/C)\eps \tag 1.32$$
where $C$ satisfies the condition (Proposition 1.12)
$$\sup_{z\in \Bbb T^2_P}\left(\int_{\Bbb T^2_P}
|\nabla g(z, \zeta)|^{p^*} d\zeta\right)^{1/p^*} \leq C\tag 1.33$$
for $p^*=(1-1/p)^{-1}<2.$
We have
$$|G^1_Dq(z)|\leq \left (\intl_{D\setminus K}+
\intl_{ K}\right )
|g^1(z,\z,D)q(\z)|d\z:=I_1(z)+I_2(z).$$
H\"older's inequality implies that
$$I_1(z)\leq \left (\intl_{D\setminus K}
|g^1(z,\z,D)|^{p*}d\z \right )^{1/p*}\times
\left (\intl_{D\setminus K}|q(\z)|^p d\z \right )^{1/p},\tag 1.34$$
and so (1.32) and (1.33) yield that
$I_1(z)<\eps$
for all $z\in D.$
For $I_2$ we have that
$$I_2(z)\leq \left (\intl_{K}
|g^1(z,\z,D)|^{p*}d\z \right )^{1/p*}\times
\left (\intl_{K}
|q(\z)|^pd\z \right )^{1/p}.\tag 1.35$$
Since $q\in L^p$ and (1.30) holds, we obtain that
$\liml_{z\ri z_0}I_2(z)=0,$
and since $\eps$ is arbitrary, these estimates
imply (1.31).
Let us note that (1.34) and (1.35) imply that
$G^1_Dq(z)$ is bounded.
The equality
$\liml_{z\ri z_0}G_Dq(z)=0 $
and boundedness of $G_Dq(z)$
can be obtained in the same way.
This and  (1.29) show that
$\liml_{z\ri z_0}q(z)=0$
and $q$ is bounded.
\qed
\enddemo

Thus (1.29) 
may be treated as
a generalization of the problem (0.15) when
$D$ is an {\sl arbitrary domain whose boundary has positive capacity}.

We recall some properties of the operator
pencil $\Cal G_D$ (see, e.g. \cite {17, Th.12.9}).

Denote by
$Spec\ D$ the {\it spectrum} of the operator pencil
$\Cal G_D$; i.e. the set of $\rho\in \BC$ for
which the operators $\Cal G_D(\r)$ have no
inverse. In particular, if $\rho \notin  Spec \ D$, every solution $q$ to
the equation (1.25) is identically zero.

It is essential for applying the cited theorem that the coefficients
$G^1_D$ and $G_D$ of the pencil be compact operators and that the
spectrum not be the whole plane.  The first assertion is Proposition 1.15; 
while the second is obvious, since the pencil is the identity
operator when $\r=0.$ From the cited theorem we obtain

\proclaim {Proposition 1.36} Let $\Cal G_D,$ $D\subset \Bbb T^2_P,$ be
the pencil (1.27).
Then $Spec\ D $
is a
discrete set (perhaps  empty) with no point
of accumulation in
any finite part of $\Bbb C.$
\ep

Here are some special properties of $Spec\ D.$ 
\proclaim {Proposition 1.37} Let $\Cal G_D,$
$D\subset \Bbb T^2_P$
be the pencil (1.27) with $Spec\ D\neq\emptyset$.
Then
\roster
\item $\{\r:\Re\r=0\}\cap Spec\ D=\varnothing$;
\item The following symmetries hold:
$$\overline {Spec\ D}:=\{\overline \r:\r\in Spec\ D\}=Spec\ D \pm
\frac {2\pi}{P}i=Spec\ D;$$
\item $Spec(D_-)=-Spec\ D$, where  $D_-$ is the  domain obtained from
$D$ by the map $z\mapsto -z$;
\item $Spec(D+z_0)=Spec\ D$ for any
$z_0\in \BT^2_P.$
\endroster

 \ep
\noindent {\it Remark.} In assertions (3) and (4) we are identifying $D$
with its image in the rectangle $R,$ extented periodically in $\Bbb C$.
\demo {Proof} We already know that
$0\not\in Spec\ D.$
If $q(z)$ is a bounded solution of the equation
$L_{\r}q(z)=0,\ z\in D$ for some $\r=it$ with $t$ real,
 then the
function $U(z)=q(\log z)|z|^{it}$ would be  harmonic and bounded in
every component $G$ of $\phi^{-1} (D)$ (see Proposition 0.10) and
would be zero quasi-everywhere on the boundary of $G$.
 Hence $U\equiv 0 $ and $q\equiv0.$
This proves assertion (2).

Similarly, given $\r \in Spec\ D$, consider the harmonic function
$U(z)=q(z) e^{\rho x}.$
Now
$$U(z)=q(z)e^{-i\frac {2\pi}{P}x}
\times e^{(\r +i\frac {2\pi}{P})x},$$
and so the function
$q_*(z):=q(z)e^{-i({2\pi}/{P})x}$, which is $P$-periodic and vanishes
q.e. on $\p D$, satisfies
$L_{\r*} q_*=0$ for $\r*=\r+ 2\pi i/{P}.$
Also, since $\overline U=
\overline q e^{\overline \r x},$
we see that  $\overline \r$ is also an eigenvalue.

 Let $D_-$ be the  domain obtained from
$D$ by the map $z\mapsto -z,$ and $\r _0,\ q_0$
be an eigenvalue and its
eigenfunction for $D_-$.
It is easy to check by changing variables that
$-\r_0,\ q(-z)$ are an eigenvalue and
eigenfunction for the domain $D.$

Finally, if $q(z)$ is an eigenfunction for an
eigenvalue $\r$ in a domain $D,$ then
$q(z+z_0)$ also satisfies (0.15) for $D:=D+z_0.$ 

\qed
\enddemo
Theorem 0.16 (to be proved in \S 6) will imply that if
$Spec\ D\neq \varnothing$, then it contains a least positive
element.  At the end of \S 6 we pose a conjecture which, if proved,
would show that it  has
considerable symmetry.

\subhead 2.\ A necessary condition that  $\Spec D\not=\varnothing$ \endsubhead

  We begin with \demo{Proof of Proposition 0.10} Let $z_0\in D$ and
$\z_0\in \phi^{-1}(z_0).$ Then $T^j\z_0 \in \phi^{-1}(z_0)$ for all
$j\in \BZ.$ Denote by $G_j,\ j\in \BZ,$ the component of
$\phi^{-1}(D)$ containing $T^j\z_0.$ Then $G_j=T^jG_0,$ because of
maximality and connectedness of $G_0$ and $G_j.$ Similarly $G_j\cap
G_l\neq \emptyset$ implies $G_j=G_l.$ We have
$$\phi^{-1}(D)=\cup_{j=-\iy}^{\iy}G_j.$$ 
If $G_j\cap G_l=\emptyset,\ j\neq l$ , then we have the 
case 1. Indeed, if there exists a curve $\hat\g\sbt D$ 
homologous to a cycle $ \g$ with $n_1\neq 0,$ then its lift $\s$ under 
$\phi^{-1}$ connects 
$\z\in G_0$ to $T^{n_1}\z\in G_{n_1},$ and, 
consequently, $G_0= G_{n_1}.$   
This is a contradiction.

Now let $\hat\g \sbt D$ be a curve corresponding to the cycle $\g$
with $\pm n_1=k\ (\geq 1).$ We can suppose that $n_1\geq 1.$ Otherwise
we replace $\hat\g$ by $-\hat\g$ with the opposite direction.  Let
$\z_0\in \phi^{-1}(\hat\g)$ and let $\s$ be the corresponding lift of
$\hat\g$ which contains $\z_0$ and $T^k\z_0.$ Denoting this component
of $G$ by $G_0,$ we have $T^kG_0=G_0.$

Let $j_{\min}\geq 1$ be the least $j\geq 1$ such that for some
$m\in\BZ$ we have $G_m\cap G_{m+j}\neq\emptyset.$ Then $j_{\min}\geq
k.$ Otherwise $G_m=G_{m+j_{\min}}$ and there exists a curve $\hat \g$
connecting $T^m\z_0$ to $T^{m+j_{\min}}\z_0.$ This means that in
$\BT^2_P$ the corresponding cycle exists with $n_1=j_{\min}<k$,
that is a contradiction. 

Since $T^kG_0=G_0,\ j_{\min}=k.$ Therefore $G_j\cap G_l=\emptyset$ for 
$0\leq j,l\leq k-1,\ l\neq j,$ i.e.,  
$G_q,\ q=0,1,...,k-1,$ are disjoint. 
Now $T^mG_0=T^qT^{lk}G_0=T^qG_0,$ and  
$T^kG_q=T^{k+q}G_0=T^qG_0=G_q, 
1\leq q\leq k-1$ .

This proves the case 2. 

\qed
\edm    
\proclaim{Proposition 2.1} Let $\Spec D\not= \varnothing$.  Then $D$
is connected on spirals.  \endproclaim

\demo {Proof} Let $q$ be an eigenfunction corresponding to an
eigenvalue $\r.$ Then $V^*(z):=q(\log z)e^{\r x}$ is a well-defined
(in general complex-valued) nontrivial harmonic function, vanishing
quasi-everywhere on the boundary of some $T$-invariant open set $G.$
The functions $\Re V^*$ and $\Im V^*$ share these properties.  We can
suppose that $V:=\Re V^*$ is positive at some point $z_0$ of a
component $G_0$ of $G.$ If $V(z)\equiv 0 $ we can replace $\Re V^*$ by
$\Im V^*.$ Otherwise $V^*(z)\equiv 0$ and hence $q(z)$ vanishes
identically. This contradicts the assumption that $q$ is an
eigenfunction.

Assume that $D$ is not connected on spirals. 
So we have  $G_i=:T^{-i}G_0\neq G_0$ for $i\neq 0,$ and 
$G_i\cap G_l=\emptyset$ for  
$i,l\in\BZ\ i\neq l$ by Proposition 0.10.

If $G_0$ is precompact in $\bC,\ V$ (and so $q$) vanishes identically
that contradicts to the assumption $V(z_0)>0.$ Thus we may assume that
each $G_i$ will have $0$ and $\infty$ in its closure (see the examples
to Proposition 0.10). Let $z_0 \in G_0,\ |z_0|=1$, and for each $j$
let $\theta_j(r)$ be the angular measure of $G_j\cap \{|z|=r\}$.
Since $G_{j}=T^{-j}G_0$, it follows that $\theta_0 (T^j
r)=\theta_j(r)$.  Let $G_j(n)=G_j\cap \{|z|=2T^{n}\}$.  By the
definition of $V$ and standard estimates on harmonic measure \cite
{20, p. 112},
$$
\aligned
V(z_0) & \leq \max_{|z|=2T^{n}} H(z)
\omega(z_0,G_0(n),G_0)\leq C T^{n\rho} \omega(z_0,G_0(n),G_0)\\
&\leq C T^{n\rho} e^{-\int_1^{T^{n}} dt/t\theta_0(t)}.
\endaligned
$$
The integral can be computed:
$$
\sum _{k=0}^{n-1}\int _{T^{k}}^{T^{k+1}}\frac{dt}{t\theta _{0} (t)}=
\sum _{k=0}^{n-1}\int _{1}^{T}\frac{dt}{t\theta_{k} (t)}=
\int _{1}^{T} \sum _{k=0}^{n-1}\frac{1}{\theta _{k} (t)}\frac{dt}{t}.
$$
Since $G_j\cap G_l=\emptyset$ for $j\neq l$ 
(Proposition 0.10), $\sum\limits_j \theta_j (r)\leq 2\pi$, so that
$n^2/ (2\pi)\leq\sum\limits_0^{n-1} 1/\theta_k(t)$,
and hence
$$
V(z_0) \leq C T^{n\rho} e^{-C_{1}n^{2}}
$$
for some constant $C_{1}>0$.
Letting $n\to\infty$, we find that $V(z_0) \leq 0$. This contradicts
the assumption  $V(z_0)>0$
and proves Proposition 2.1.
\qed
\enddemo

\subhead 3.\ Proof of Theorem 0.5\endsubhead
 Our approach
uses the calculus of positive harmonic functions introduced by
R. S.~Martin \cite {18} and popularized in the thesis of the late
B.~Kjellberg \cite {13}.  Martin's insight was to consider limits of
ratios of the type
$$
h_n(z)={\omega(z,G\cap E_n, G)\over \omega(z_0,G\cap E_n, G) }\tag 3.1
$$
where $z_0$ is fixed in $G\cap \{|z|=1\}$ and the $\{E_n\}$ are sets
of positive capacity tending to $\infty$.  We may assume, as is
customary, that any function $v(z)\in \sP$ is zero for $z$ not in $G$.
Let $\Delta_\infty$ be the cluster set (in the topology of uniform
convergence on compact subsets of $G$) of the functions $h_n$ which is
obtained by letting the $\{E_n\}$ tend to infinity.  It consists of
positive harmonic functions on $G$ which are $1$ at $z_0$.  Notice
that in general $\Delta_\infty$ need not be contained in $\sP$,
because a limit function need not vanish q.e. on $\partial G$.

Let us introduce some examples, where again $G$ and $D$ are related by
(0.9).  Let $G_\rho$ be the sector $\{|\arg z| < {\pi}/{2\rho}\}$.
Then the cone $\sP$ of Theorem 0.5 consists of positive multiples of
the function $u_\rho(r,\theta)=r^\rho \cos(\rho\theta)$, so that $\sP$
has dimension one and $\sP=\sF$.

An illuminating example of a 2--homogeneous set  
for which $\Cal F\sbt\Cal P$ 
is the set $\Omega_0$:
$$ \Omega_0=\{y > 0\}\backslash\bigcup^\infty_{n=-\infty} \{y=2^n,\
-\infty < x < 0\}, $$ the upper half--plane with a sequence of
horizontal rays deleted (this example provided the original motivation
for this section).  If $E_n\to\infty$ inside the first quadrant, the
family associated to $\{E_n\}$ by (3.1) will converge to a function
$u\in \Delta_\infty$ which is also in $\sF$, and Theorem 0.5 implies
that the positive multiples of $u$ span all of $\sF$.  However, if the
$\{E_n\}$ tend to $-\infty$ through one of the horizontal channels,
say ${\Cal C}_k= \{x<0, 2^k<y<2^{k+1}\}$, then the $h_n$ converge to a
function $u_k$ in $\Delta_\infty$, which will also be in $\sP$, but
which has infinite order and hence is not in $\sF$.  However, $u_k$
will be bounded outside the given channel ${\Cal C}_k$, and so, if
$u_j$ and $u_k$ are associated to two distinct channels, they will be
linearly independent. So in this case $\sP$ contains at least
countably many linearly independent functions.

The remainder of this section is devoted to the proof of Theorem 0.5 .

Let us note that assertion (2) is an easy corollary of (1).
Indeed, the function
 $H(Tz)$ is in $\sF$ along with $H(z),$ and hence $H(Tz)=
cH(z)$ with
$$c:=H(Tz_0)/H(z_0)=H(Tz)/H(z) $$ Hence $c$ does not depend on $z_0$.
It does not change if we replace $H$ for any $v\in \sF$ because $v=aH$
with a constant $a .$ Define $\r(G)$ by $c=T^{\rho(G)},$ and note that
$v$ satisfies (0.7).  It is clear that $\r (G)$ cannot be negative or
zero because in such case $H(z)$ would be bounded and hence vanish
identically.

Let us also remark that we can replace $G$ with
$re^{i\psi}G$ for
any $re^{i\psi}$ without loss of generality. Indeed,
$G'=re^{i\psi}G$ is also $T$ homogeneous and equality
$v_1(\z):=v(\z/re^{i\psi}), \z\in G'$ generates
one-to-one maps $\sP_{G}\mapsto\sP_{G'}$ and
$\sF_G\mapsto\sF_{G'},$ and $\r(G')=\r(G).$

Next we produce $H(z)$ and prove uniqueness.
Our first result, Lemma 3.3,
gives a concrete way to characterize $\sF$ in $\sP$.
 We set
$
T_n=\{|z|=T^n\},
$
and suppose  as in Theorem 0.5
that $z_0\in G,\ |z_0|=1$, has been fixed.

 For $v\in\sP$, set
$$
\beta(v)=\limsup_{n\to\infty} \{\max_{\zeta\in T_n}v(\zeta)
\omega(z_0,T_n)\}. \tag3.2 $$

\proclaim{Lemma 3.3} Let $v\in\sP$ and $\beta(v)=\infty$.
Then $v$ has infinite order.
\endproclaim

\demo{Proof} Partition $T_n \cap G$ into $I_{n}=I_n (\varepsilon )$ and
$J_{n}=J_n (\varepsilon )$, where $J_n=T_n\cap \{z\in G: d(z,\partial
G) < \varepsilon T^n\}$ and $I_n=(T_n\cap G)\setminus J_n$.
In the sequel, we will often omit the dependency of $I_{n}$ and
$J_{n}$ on $\varepsilon$.

We need the following fact:
$$
\lim_{\var\to 0}\sup_{z\in T_0} \omega(z,J_{1} (\varepsilon))=0.\tag3.4
$$
Although this is easy to verify for most domains, in the generality in
which we are working we need a more careful justification.

\demo{Proof of (3.4)}Let us replace $G$ by $rG$
with $r$ chosen so that
$$\om(z_0,\partial G\cap T_{1}) =0,\tag 3.5 $$
with $z_0$ a base point in $G\cap T_0$.
This is possible because the function
$f (t)=\om(z_0,\partial G\cap \{1<|z|<t \})$ is monotone increasing
with $t$ and hence has a  
dense set of continuity points in any interval.
Notice that this does
not restrict generality because of the remarks made above at the
beginning of proof Theorem 0.5.

Set  $h_\varepsilon(z):=\om (z,J_1(\varepsilon),G).$
Since $\{h_\varepsilon,\ \varepsilon >0\}$ is a
bounded family of harmonic functions it is a
{\it normal} family. Hence for an arbitrary
sequence $\ve_j\ri
 0$ there exists a subsequence
$\{\ve_{j'}\}$ such that $h_{\ve_{j'}}$
 converges
uniformly on every compact subset of $G$ to
a harmonic function $h(z).$

Suppose   for some
sequence $\{\ve_{j}\}$ there exists a subsequence
such that $h(z)\not\equiv 0$ and hence
$$h(z_0)> 0\tag 3.6$$
by minimum principle.

For every $h_{\ve}$ we have the inequality $h_{\ve}(z)\leq \om
(z,T_1,G), z\in G.$ Thus $h$ itself satisfies the same inequality and
hence $\lim \limits_{z\ri\z}h(z)=0$ at any regular point $\z\in
\partial G$, except, possibly, points of $E:=\partial G\cap T_{1};$
i.e. q.e. on $\p G\setminus E.$ Denote as $\partial_I G$ the the set
of irregular points in $\p G$. Since $\om (\p_IG)=0$, we combine these
estimates with the maximum principle and (3.5) to deduce that
$$h(z_0)\leq \max\{\z\in \p G\setminus
(E\cup \p_I G):\limsup\limits_{z\ri\z}h(z)\}\cdot
\om (z_0,\p G\setminus (E\cup \p_I G))$$
$$+1\cdot\om (z_0,\p_I G)+1\cdot\om (z_0,E)= 0 + 1\cdot 0 + 1\cdot 0,$$
and this contradicts
(3.6).
Thus $h(z)\equiv 0$ and (3.4) is proved.
\qed\edm

With $z_0\in G\cap \{|z|=1\}$ fixed as above, 
choose a path $\gamma$
joining $z_{0}$ to $Tz_{0}$ in $G,$ and let $\Omega$ be open with
compact closure in $G$ such that $\Omega\supset \gamma$. We may then
take $\tau$ to be the Harnack constant $\tau(\gamma, \Omega).$ So
if $u$ is positive and harmonic in $G\cap \{|z|<R \}$
with $R>T$ large enough so that $\Omega\subset\{|z|<R \}$, then
$$
\tau^{-1} u(z_0)\leq u(Tz_0)\leq \tau u(z_0).\tag3.7
$$

We now
choose $\var_0>0$ so small in the definitions of $I_n$ and $J_n$
so that $T^n z_0\in I_{n}$ and, using (3.4), so that
$$
\sup_{z\in T_0} \omega(z,J_1)\leq {1\over 2\tau}.\tag3.8
$$

In order to appreciate the significance of (3.2), we show that if
$v\in\sP$, then
$$
v(T^n z_0)\omega (z_0,T_n) < Bv(z_0),\tag3.9
$$
for some constant $B= B(G,z_{0})<\infty.$ 
Thus the condition $\beta(v) = \infty$, forces  $v$ to grow rapidly away from
the orbit of $z_0$.

To show (3.9) we first note that,
$$
\aligned
\omega(z_0,J_n)&\leq \omega(z_0,T_{n-1})\sup_{\zeta\in T_{n-1}}
\omega(\zeta,J_n)\\
&\leq \tau\omega(z_0,T_n)\sup_{\zeta\in T_0}\omega(\zeta,J_1)\leq
{1\over 2}\omega(z_0,T_n).
\endaligned
\tag3.10
$$
where the first inequality follows from the strong Markov property;
the second one uses (3.7) and $T$--homogeneity; and the last one uses (3.8).
We remark that (3.10) holds  only for $n$ large, i.e. for
$n\geq n_0$ where $n_0$ depends on $G$ and $z_{0}$ so that (3.7) can be used.

Since $\var_0 > 0$ has been fixed, it follows from  Harnack's inequality on
$I_{n}=I_{n} (\varepsilon _{0})$
and $T$-homogeneity that
there exists a constant  $0 <b_0=b(\var_0,G)<1$ with
$\min\limits_{I_n} v(z)\geq b_0
v(T^nz_{0})$ for all $n$.
Thus, we deduce for $n\geq n_0$ that
$$
\aligned
v(z_0)&\geq \omega(z_0,I_n)\min_{I_n} v(\zeta)\\
&\geq {1\over 2}\omega (z_0,T_n)\min_{I_n} v(\zeta)\\
&\geq {b_0\over 2}\omega (z_0,T_n) v(T^n z_0),
\endaligned\tag3.11
$$
where the first inequality follows by the maximum principle on the
region $G\setminus I_{n}$; the second uses (3.10) and
subadditivity of harmonic measures, i.e. the fact, which follows from
the maximum principle, that
$$
\omega (z_{0},T_{n},G\setminus T_{n})\leq \omega
(z_{0},I_{n},G\setminus I_{n})+ \omega (z_{0},J_{n},G\setminus
J_{n});
$$
and the last inequality in (3.11) uses Harnack's
Inequality. Thus (3.9) is proved.

We can now finish the proof of Lemma 3.3.
Take $z\in G,\ |z| \leq T^n,$ and set $M_n=\max\limits_{T_n} v(\zeta)$.
By assumption we have
$$
\limsup M_n \omega(z_0,T_n)=\infty.
$$

Given $S>1$,
(3.4) again implies that we may decrease $\var$ in
the definitions of $I_n$ and $J_n$ so that
$$
\sup_{T_0} \omega(z,J_1) < {1\over 2\tau S},
$$
and so that (3.10) and (3.11) still hold, with a different
constant $0<b=b(\var,G)<1$ instead of $b_0$. Then,
for $|z| < T^n,\ z\in G,$
$$
\aligned
v(z)&\leq M_{n+1}\omega (z, J_{n+1})+\max_{I_{n+1}} v(\zeta)\\
&\leq M_{n+1} {1\over 2\tau S} +\max_{I_{n+1}} v(\zeta)\\
&\leq M_{n+1} {1\over 2\tau S}+(\tau/b)v (T^n z_0)\\
&\leq M_{n+1} {1\over 2\tau S}+(2\tau/b^2)v(z_0)(\omega (z_0,T_n))^{-1}
\endaligned
$$
where  the first inequality follows from the maximum principle on
$G\cap \{|z|<T^{n} \}$; the second one uses our choice of $\varepsilon
$; the third one uses Harnack's inequality on $I_{n+1}$ (as already done just
before (3.11)), as well as (3.7); and the final
inequality follows (3.11).
Taking the supremum over all $z$'s for which the above inequality 
holds, multiplying both sides by $\omega
(z_{0},T_{n})$, and using (3.7) on $\om (z,T_{n})$, we obtain that
$$
 M_n \omega(z_0,T_n)\leq M_{n+1} \omega(z_0,T_{n+1}) {1\over 2S}+Av(z_0),
$$
with  $A=A(\var)$,
and so by iterating this inequality for $k=1, 2, 3 \dots$,
$$
 M_n \omega(z_0,T_n)\leq {1\over (2S)^k} M_{n+k}\omega
(z_0,T_{n+k})+Av(z_0)\sum_0^{k-1} {1\over (2S)^j}.
\tag3.12
$$
Since $\beta(v)=\infty$, we may
choose $n=n_1$ so large that
$$
M_{n_1} \omega(z_0,T_{n_1}) > \frac{4SAv(z_0)}{2S-1}>2Av(z_0),
$$
and then (3.12) implies for each $k > 1$ that
$$
M_{n_1+k} \omega(z_0,T_{n_1+k})\geq (2S)^k\left(M_{n_1}\omega(z_0,T_{n_1})
-\frac{Av(z_0)}{1-(1/(2S))}\right)\geq Av(z_0)(2S)^k.
$$
Since $S$ is arbitrary, $v$ must have infinite order.
Lemma 3.3 is proved.
\qed
\enddemo

We now construct functions in $\Cal F$.
Choose $E\subset G\cap T_0$ with $|E| > 0$ and (cf.~(3.1)) let $\sH$
consist of all normal limits of the family of functions
$$
h_n(z)={\omega(z,T^n E)\over \omega(z_0,T^n E) }\qquad
(z\in G\cap \{|z| \leq T^n\}\quad n=1,2,3,\dotsc ).\tag3.13
$$

\proclaim{Lemma 3.14} Let $\sH$ be as above.
Then $\sH\subset\sF$.
\endproclaim

\demo{Proof} We first show that for $m$ large,
$$
\sup_{T_0} \omega(\zeta,T^m E)\leq D\omega(z_0,T^m E),\tag3.15
$$
for a constant $D>0$ which only depends on the domain $G$.
Recall that in the proof of Lemma 3.3 we created $\{I_n=I_{n}
(\varepsilon _{0})\},\ \{J_n=J_{n} (\varepsilon _{0})\}$ so that (3.8)
holds, and recall the constants $\tau$ and $b_0$ as well.
If $S_m=\sup_{\zeta \in J_{0}} \omega(\zeta,T^m E)$, then
$$
\aligned
S_m&\leq \sup_{I_1} \omega(z,T^m E)+{1\over 2\tau}\sup_{J_1}\omega
(\zeta, T^m E)\\
&\leq(\tau/b_0)\omega(z_0,T^m E)+{1\over 2\tau} S_{m-1}\\
&\leq(\tau/b_0)\omega (z_0,T^m E)+{1\over 2} S_{m-1} {\omega(z_0,T^m E)\over
\omega(z_0,T^{m-1} E)}
\endaligned\tag3.16
$$
where the first inequality follows from the maximum principle on
$G\cap \{|z|<T \}$ and (3.8); the second one uses Harnack's inequality
on $I_{1}$ followed by (3.7) and the definition of $S_{m-1}$; and  the
last line follows from (3.7) and homogeneity.

We deduce  from (3.16) that
$$
\aligned
{S_m\over \omega(z_0,T^m E)}&\leq \tau/b_0+{1\over 2} {S_{m-1}\over
\omega(z_0,T^{m-1} E)}\leq \ldots\\
&\leq (\tau/b_0) \sum_0^{m-m_1-1} 2^{-j}+2^{-(m-m_1)} {S_{m_1}\over
\omega(z_0,T^{m_1} E)}\leq 4\tau/b_0,
\endaligned
$$
for $m$ large enough.  Harnack's inequality on $I_{0}$ yields that
$$\sup_{\zeta\in I_0}\omega(\zeta,T^mE)\leq(1/b_0)\omega(z_0,T^mE),$$  and so (3.15)
holds with $D=4\tau/b_0$.

Consider now the functions $\{h_n\}$ of (3.13), and let $|z|\leq T^n$.
Then for $m>n$ ($m$ much larger than $n$),
$$
\aligned
\omega(z,T^m E) & \leq \omega(z,T_n)\sup_{\zeta\in T_n}\omega(\zeta,T^mE)\\
&\leq D\omega (z,T_n) \omega (T^n z_0, T^m E)\\
&\leq D\omega(z,T_n) \tau^n \omega (z_0,T^m E),
\endaligned
$$
where the first inequality follows from the maximum principle on
$G\cap \{|z|<T^{n} \}$; the second one uses (3.15) and homogeneity;
and the last one (3.7) $n$ times.
Hence  any normal limit $h$ of the $h_n$ must satisfy
$$
h(z)\leq D\tau^n \omega (z,T_n),
$$
for $|z|\leq T^{n}$, and so  $h$
is locally bounded and vanishes q.e. near
each finite boundary point of $G$.
If we set $\tau = T^\ell$, then with the notations of (0.6),
$M(T^n,h) \leq C \tau^n = C(T^n)^\ell$,
so that   $h\in \sF$ (hence by Lemma 3.3, $\beta (h) < \infty$).
Thus $\Cal F \neq \varnothing$.
\qed
\enddemo
Finally, we show that $\sF$ is one--dimensional.
Let $H$ be any limit function of the family (3.1) and let  $\langle H
\rangle$ consist of  all positive multiples of $H$.
Clearly $\langle H \rangle \subset\sF$.
\proclaim{Theorem 3.17} $\langle H\rangle =\sF$.
\endproclaim

\remark{Remark} If $G\cap \{1 < |z| < T\}$ were a Lipschitz domain,
this would be a consequence of the boundary Harnack principle (cf. \S5).
What follows is a replacement for this principle.
\endremark

The strategy to prove uniqueness is as follows: we first construct
some auxiliary functions $V_\var\in \sF$,
namely for each $\var>0$ small enough, we use
the  partition $\{I_n (\varepsilon ),J_n (\varepsilon )\}$
that was discussed in the
proof of Lemma 3.3., and  then we produce $V_\var$ as a certain
sublimit of the ratios (3.13) with $E=I_0$; once this is done we show
that every function in $\sF$ is comparable to the functions $V_\var$,
and then we conclude using a standard argument.
First we need a lemma.

\proclaim{Lemma 3.18} There exists $\var_{1}=\var_1(G) > 0$ and $K
=K(G)> 1$ so that,
whenever $\{I_n=I_{n} (\varepsilon )\}$ and $\{J_n=J_{n} (\varepsilon
)\}$ are created as in the proof of Lemma
3.3 for any fixed  $0<\var<\var_{1}$, and $U$ is defined as
$$
U(z)=\lim_{k\to\infty} {\omega(z,J_{m_k})\over \omega(z_0,J_{m_k})}
$$
$($where $m_1 < m_2 < \ldots$ $)$,  there is a limit function $V$ of
the corresponding family
$\{V_{m_{k}}\}=\{\omega(z,I_{m_k})/\omega(z_0,I_{m_k})\}$
with
$$
U(z) < KV(z)
$$
for all $z\in G$.
\endproclaim

\demo{Proof} With $\tau, B, D$ from (3.7), (3.9) and
(3.15) choose $\var_{1}$ so that for $\var <\var_{1}$
$$
\sup_{T_0}\omega(z,J_1) < {1\over 2BD\tau}.
$$
Since $\tau, B$, and $D$ only depend on the domain $G$, also
$\var_1$ only depends on $G$.
Fix $m > n$ and let $|z| < T^{n-1}$.
We shall analyze the inequality
$$
\omega(z,J_m)\leq
\omega(z,I_n)\sup_{I_n}\omega(\zeta,J_m)+\omega(z,J_n)\sup_{J_n}
\omega(\zeta,J_m),\tag3.19
$$
which follows from the maximum principle on $G\cap \{|z|<T^{n} \}$.

Consider the first term on the right.
With $M=\min\limits_{I_1}
\omega(\zeta,I_{0}),\  0 < M < 1$, we have
$$
\omega(z,I_n)\leq M^{-1} \omega(z,I_{n-1}),
$$
by the maximum principle on $G\setminus (I_{n-1}\cup I_{n})$ and
homogeneity.

By Harnack's inequality on $I_n\cup {T^{n-1}z_{0}}$,
$$
\sup_{I_n}\omega(\zeta ,J_m) < A\omega (T^{n-1} z_0, J_m)\qquad (A=A(\var)),
$$
so that the first term on the right side of (3.19) is at most
$$A M^{-1} \omega(z, I_{n-1})\omega (T^{n-1} z_0, J_{m}).$$
As for the second term,
$$
\aligned
& \omega(z,J_n)\leq \omega(z, I_{n-1})\sup_{I_{n-1}} \omega(\zeta, J_n)\\
&\quad +\omega(z, J_{n-1}) \sup_{J_{n-1}} \omega (w, J_n)\\
&\quad \leq \omega (z, I_{n-1})+{1\over 2BD\tau}\omega (z, J_{n-1})
\endaligned
$$
where the first inequality follows from the maximum principle on
$G\cap \{|z|<T^{n-1} \}$; and the second one our choice of
$\varepsilon $.
On the other hand,
by (3.15), homogeneity,  and (3.7), we have
$$
\sup_{J_n} \omega(\zeta, J_m)\leq D\tau \omega (T^{n-1} z_0, J_m),
$$
and so  we deduce that the second term is at most
$$
\omega(z,I_{n-1}) D\tau\omega (T^{n-1} z_0, J_m)+{1\over 2BD\tau}
\omega(z,J_{n-1})D\tau \omega(T^{n-1} z_0, J_m).
$$
Combining these estimates, we obtain that
$$
\aligned
\omega(z,J_m)&\leq (M^{-1} A+D\tau) \omega(z,I_{n-1})\omega (T^{n-1}
z_0, J_m)\\
&+ \omega(z, J_{n-1})\cdot {1\over 2B}\omega (T^{n-1} z_0,J_m).
\endaligned
$$
Divide both sides by $\omega(z_0,J_m)$, and let $m\to\infty$ appropriately.
It follows that
$$
U(z)\leq (M^{-1} A+D\tau)\omega (z, I_{n-1}) U(T^{n-1} z_0)+{1\over
2B} \omega(z,J_{n-1}) U(T^{n-1} z_0).
$$
Note that by Lemma 3.14, $U$ belongs to $\Cal P$,
and (3.9) applied to $U$ refines the last estimate to
$$
U(z)\leq (M^{-1} A+D\tau)B { \omega(z,I_{n-1})\over \omega (z_0,
I_{n-1})}+{1\over 2} {\omega(z,J_{n-1})\over \omega(z_0, J_{n-1})},
$$
where we used the fact that $U (z_{0})=1$.
Now let $n-1$ tend to infinity along an appropriate subsequence of the
$\{m_k\}$.
Then
$$
U(z)\leq (M^{-1} A+D\tau)B V(z)+{1\over 2} U(z),
$$
and so we may take $K=2 (M^{-1} A+D\tau)B$. Lemma 3.18 is proved.\qed\enddemo

Suppose  that for every $0<\var<\var_1$ we construct the partition
$\{I_n=I_{n} (\varepsilon ),J_n=J_{n} (\varepsilon )\}.$
Normal families then produce a function
$U_{\varepsilon }$, and hence also a function $V_{\varepsilon}$, as in
Lemma 3.18.

We next study an expression complementary to $\beta (v)$ in (3.2).

\proclaim{Lemma 3.20} Let $v\in\sF$.
Then there exists $\var_v > 0,$ so that if $\var < \var_v$ and
$\{I_n=I_{n} (\varepsilon )\},\ \{J_n=J_{n} (\varepsilon )\}$
are constructed as in Lemma 3.3, there is a constant $A_\varepsilon>1$
such that
$$
\omega(z_0,I_n)\min_{I_n} v(\zeta) > v(z_0)/A_{\var}\qquad (n=0,1,2,\dots ).
$$
\endproclaim

\demo{Proof} Since $\beta(v) < \infty$,
$$
S(v) = \sup_n\max_{T_n} v(\zeta) \omega(z_0, T_n) < \infty.
\tag3.21
$$
By (3.4), we may choose $\var_v$ so that if $\var < \var_v$, then
$$
\sup_{T_0}\omega(z,J_1) < {v(z_0)\over 2\tau S(v)}.
$$
By the argument which yielded (3.10), we have
$$
\omega(z_{0}, J_n)\leq {1\over 2S(v)} \omega(z_0, T_n).
$$
Thus,
$$\aligned
v(z_0)&\leq \omega(z_0, J_n)\max_{J_n} v(\zeta)+\omega (z_0,
I_n)\max_{I_n} v(\zeta)\\
&\leq {v(z_0)\over 2S(v)} \omega(z_{0},T_n)\max_{T_n} v(\zeta)+A_\var\omega
(z_0,I_n)\min_{I_n} v(\zeta)\\
&\leq {v(z_0)\over 2}+A_\var\omega (z_0,I_n)\min_{I_{n}} v(\zeta),
\endaligned
$$
where the first inequality follows from
the maximum principle on $G\cap \{|z|<T^{n} \}$; the second one uses
our choice of $\varepsilon $ and Harnack's inequality on $I_{n}$; and
the last one uses (3.21).
Lemma 3.20 is proved provided $A_{\varepsilon }$ is changed to
$2A_{\varepsilon }$.
\qed
\enddemo

Now we show
that
the class $\sF\subset\sP$ is one-dimensional; this will
prove Theorem 3.17  and thus Theorem 0.5.
By choosing $0<\var<\min\{\var_1,\var_v\}$
and letting $V_\var$ be the function constructed after Lemma
3.18, we show that any $v \in \sF$ satisfies
$$
C^{-1} v(z)\leq V_{\varepsilon }(z)\leq Cv(z)\tag3.22
$$
for all $z\in G$.
Once (3.22) is established, we use a
now--standard argument of Kjellberg, which we recall at the end of the
proof for completeness. On ``nice'' domains (cf. \S 5) equation (3.22)
would follow automatically from the boundary
Harnack principle and homogeneity.

By the maximum principle, for all $z\in G,\ |z|\leq T^n$,
$$
v(z)\geq {\omega(z,I_n)\over \omega(z_0,I_n)}
\omega(z_0,I_n)\min_{I_n} v(\zeta),
$$
so, by Lemma 3.20 and the construction of $V_\var$, $v(z)\geq
(v(z_0)/A_\var) V_\var(z)$.
So the second inequality of (3.22) is proved.

On the other hand, if $|z| < T^n$,
$$
\aligned
v(z)&\leq {\omega(z,I_n)\over \omega(z_0,I_n)}\max_{I_n}
v(\zeta)\omega(z_0,I_n) + {\omega(z,J_n)\over \omega(z_0,J_n)}
\max_{J_n} v(\zeta)\omega(z_0,J_n)\\
&\leq \left[{\omega(z,I_n)\over \omega(z_0,I_n)}+{\omega(z,J_n)\over
\omega(z_0,J_n)}\right]\max_{T_n} v(\zeta)\omega (z_0,T_n),
\endaligned
$$
where both inequalities follow from the maximum principle on
$G\cap \{|z|<T^{n} \}$.
Since $v\in \sF$, Lemma 3.3 implies that
$\limsup_{n\to\infty}\max_{T_n} v(\zeta)\omega
(z_0,T_n)=\beta (v)<\infty$. Hence, by Lemma 3.18,
$v(z)\leq CV_\var(z)$, and (3.22) follows.

We now repeat Kjellberg's argument.
Suppose $v_1,v_2 \in \sF$. Two
applications of (3.22) imply that $v_1$ and $v_2$ i
must also be comparable. Set
$$
m=\inf_{G} {v_2(z)\over v_1(z)}\in[C^{-1},C].
$$
We claim that $v_2-m v_1\equiv 0$.
Assume the contrary.
Then $w=v_2-m v_1\in \sF$, and by (3.22)
$w(z)/v_1(z) > 1/C_{1}$ for every  $z\in G$, for
some other constant  $C_{1}>1$.
But
$$
\inf_{z\in G}\frac{w(z)}{v_1(z)} = m-m=0.
$$
This is a contradiction.
Therefore $w=0$ and $v_2 \equiv mv_1$. So $\sF$ consists of positive
multiples of a single function.

\subhead 4.\ $\r (G)$ as an order of growth and decay of
harmonic functions in $G$ \endsubhead
Let $H$  be the unique function in $\Cal F$
as determined by Theorem 0.5.
Set
$$
M(r,H)=\sup_\theta\{H(re^{i\theta}): z=re^{i\theta}\in D\}.$$

\proclaim{Proposition 4.1} The function $H$ is related to $\rho (G)$ by
$$
\rho(G) = \lim_{r\to \infty}\frac{\log M(r,H)}{\log r} =\lim_{r\to 0}\frac{\log M(r,H)}{\log r},\tag4.2
$$
where the limits exist and are positive.
\endproclaim

\demo{Proof} From  (0.7) we obtain that
$M(Tr,H)=T^{\r (G)}M(r,H)$, which implies
$$\log M(T^nr,H)=n\r (G)\log T +\log M(r,H).$$
This means that
$$\r (G)=\lim_{n\to \infty}\frac{\log M(T^n,H)} {n\log T} =\lim_{n\to
-\iy}\frac{\log M(T^n,H)} {n\log T}.$$ Note that $M(r,H)$ is
increasing with $r$, by the maximum principle. This and monotonicity
of $\log r$ imply (4.2).  From Theorem 0.5 $\r (G)>0.$ \qed \enddemo
The constant $\rho(G)$ is also intimately related to the decay of
harmonic measure. Recall that $\om (z,T_n)$ is the harmonic measure of
$T_n:=\{z:|z|=T^n\}$ with respect to $G.$

\proclaim{Proposition 4.3}For each $z\in G$
there corresponds $C=C(z) > 1$ such that
$$
(1/C) T^{-n\r (G)}\leq \omega(z,T_n)\leq
C T^{-n\r (G)}\quad (n\geq 1).
\tag 4.4$$

In particular, for each $z\in G$,
$$
\lim_{r\to\infty}\frac{\log \omega(z,\{|\zeta|=r\})}{\log(1/r)}=
\lim_{n\to\infty}\frac{\log\omega(z,T_n)}{\log T^{-n}}=\rho (G)\tag 4.5$$
\ep
\demo{Proof} Let $H\in\sF$ be from Theorem 0.5 with $H(z_0)=1$.
Then (3.9) yields that
$$\omega(z_0,T_n)\leq B[H (T^n z_0)]^{-1}=
B T^{-n\r (G)},
$$
which, using Harnack's inequality to compare $\omega (z,T_{n})$ to
$\omega (z_{0},T_{n})$, 
yields the right-hand estimate of (4.4), while  Lemma 3.20 shows
that
$$
\omega(z_0, T_n) H(T^n z_0)\geq
\omega(z_0,I_n)\min_{I_n} H(\zeta)\geq A,
$$
leading to the reverse inequality. So (4.4) is
proved.

The second equality in (4.5) follows obviously
from (4.4). The first one follows from the second one and
monotonicity of $\om(z,\{|\zeta|=r\})$ and
$\log(1/r)$ in $r.$
\qed
\enddemo

\subheading {4.6} The set-function $\r(G)$ has a classical
interpretation in the situation that $G$ is simply-connected and
$T$-invariant.  Let $I_0$ be the interval of $G\cap T_0$ that
separates zero from infinity in $G$ and set $I_1=T I_0.$ Since $G$ is
simply-connected, $I_0$ divides $G$ into two components, as does
$I_1.$ Let $R(I_0, I_1)$ be the quadrilateral which is the component
of $G$ having $I_0$ and $I_1$ on its boundary, and let $\mod
R(I_0,I_1)$ be its conformal modulus. Namely $\mod R(I_0,I_1)$ is the
length $L$ of the unique rectangle $Q=\{z=x+iy:0<x<L, 0<y<1 \}$, which
can be obtained by mapping $R (I_{0},I_{1})$ to $Q$ conformally in
such a way that $I_{0}$ and $I_{1}$ are mapped to the two vertical
sides respectively.

Let us make the following construction.
Denote  $G_0:=R(I_0,I_1)$ and set 
$G_n:=T^nR(I_0,I_1)=R(I_n,I_{n+1}),\ I_n=T^nI_0\ n\in \BZ.$
Replace every $G_{2k-1},\ k\in \BZ$ by 
$G^*_{2k},$ the quadrilateral that is symmetric to $G_{2k}$ with
respect to the arc $I_{2k}.$ We
obtain a new $T^2$-homogeneous domain $G^S$,
or two-sheeted Riemann surface, because 
a point can be covered by some  $G_{2k}$ and possibly by some  $G^*_{2m}.$

We do not develop our theory for Riemann surfaces, so we confine
ourselves to domains $G$ with a {\it separating circle\,}: there
exists a circle that intersects $G$ on one arc. Without loss of
generality take that to be $I_0.$ Then $G_n$ and $G^*_{n-1}$ have the
only common arc $I_n=T^nI_0 \in G.$ So if $G$ has a separating circle
$G^S$ is a plane 
domain and vice versa.  
\proclaim { Theorem 4.7} If G is a
domain with a separating circle
$$\r(G^S)=(\pi/P)\mod R(I_0, I_1).\tag4.8$$ \ep \demo {Proof} Let
$w=f(z)$ map $R:=G_0$ conformally to the rectangle $(0,c)\times (0,
\pi)$ (for a unique $c$) so that $I_0$ and $I_1$ correspond
respectively to the vertical sides.

  By reflection, we may
extend $f$ to map $R^*=G^*_{-1}$ to the rectangle $(-c,0)\times (0,
\pi),$ and in the same way extend $f$ to map $G^S$ to the strip
$\{0<v<\pi\}$ so that

$$ f(T^2z)=2c + f(z),\ \Re f(T^2z) =2c + \Re f(z),$$
On iterating, we have that  $\Re f(T^{2n}z) =2cn + \Re f(z).$

 The function 
$H(z):=\Im e^{f(z)}=
e^{\Re f(z)}\sin\Im f(z)$
is positive and harmonic within $G^S$ and equal 
to zero at every regular point of boundary. One can check that it has
a finite order. By (4.2)  
$$\rho (G^S) = \lim_{n\to\infty} \frac{2cn + O(1)}{2n \log T}= \lim_{n\to
\infty} \frac{2cn + O(1)}{2nP}=\frac {c}{P}=
\frac{\pi}{P}\cdot\mod R(I_0,I_1).$$

\qed
\edm
\proclaim {Corollary 4.9} If $G$ is a T-homogeneous domain with a
separating circle
which is also symmetric with respect to
reflection in this circle, then
$$\r(G)=(\pi/P)\cdot\mod R(I_0, I_1).$$ \ep

The set-function $\r(G)$ has  several other
interpretations in the situation that $G$ is simply-connected and
$T$-invariant (see Lemma 6.4 and Theorem 6.6 of \cite {19}).  Thus
 Proposition 4.3
can also be formulated in terms of
the
growth of the distance in the hyperbolic metric of $G$ between $T_{0}$ and
$T_{n}$, as well as the growth of the extremal distance  between $T_{0}$ and
$T_{n}$. These latter  objects make sense also for the non-simply-connected
domains $G$ (because we assume the
capacity of the boundary to be non-zero), although explicit computation
may be more difficult.

Here we discuss one such result related to
\cite {19}.

\proclaim{Proposition 4.10}
Let $G$ be simply connected and $T$-invariant for some $T>1$. Let
$I_{0}$ be an arc of $G\cap \{|z|=1 \}$ which
separates zero from $\iy$ in $G$  and let $I_{n}=T^{n}I_{0}$.
Then
$$\rho (G) = \lim_{n\to\infty}\frac{\pi}{\log T} \frac{d_{G}(I_0,I_n)}{n}.$$
where $d_{G} (I_{0},I_{n})$ is the extremal distance between the two
crosscuts $I_{0}$ and $I_{n}$ in the quadrilateral formed by them in $G$.
\endproclaim
The same formula holds if one uses all the arcs $T_{0}=G\cap \{|z|=1 \}$
and $T_{n}=T^{n}T_{0}$, however the proof is much more delicate and
can be deduced from the proof of Claim 6.7 of \cite {19}.  The problem is that
in general, with the notations of the proof of Proposition 4.10 below,
$\psi^{-1}(T_0)$ can be quite pathological, so that both $0$ and
$\infty$ are in its closure.
\demo{Proof of Proposition 4.10}
As constructed in Lemma 6.4 of \cite {19} let $\psi$ be a conformal map of
the upper half-plane $\bH$ onto $G$ such that
$$
\psi (tz)=T\psi (z)
$$
for some $t>1$. Let $J_{0}=\psi^{-1} (I_{0})$ and $J_{n}=\psi^{-1}
(I_{n})=t^{n}J_{0}$. By properties of conformal maps, $J_{0}$ is a
Jordan arc in $\bH$ with two end-points $a,b \in \bR\setminus\{0 \}$.
Therefore $M=\max\{|z|:\ z\in J_{0} \}$ and $m=\min \{|z|:\ z\in J_{0}
\}$ are well defined with $0<m\leq M<\infty$. Let $C_{r}$ be the arc
$\{|z|=r \}\cap \bH$.
By conformal invariance, $d_G(I_0, I_n) = d_{\bH}(J_0,J_n),$ so
for $n$ sufficiently large,
$$
\frac{1}{\pi} (n\log t-\log \frac{M}{m})
\leq d_{G} (I_{0},I_{n})\leq \frac{1}{\pi} (n\log t+\log \frac{M}{m}).
$$
Hence
$$
\lim_{n\to\infty}\pi\frac{d_{G}(I_0,I_n)}{n}=\log t.
$$
The function $H$ of Theorem 0.5 is a positive
multiple of $h (z)=\Im\psi^{-1} (z),$ and $h (Tz)=th (z)$. Thus
(0.7) yields that
$$
\rho (G) =\frac{\log t}{\log T}.
$$
\qed
\enddemo

\subhead 5. A sufficient condition that $\sF=\sP$\endsubhead We saw at
the beginning of \S3 that in general $\sF\subset\sP$, with strict
inclusion possible.  Suppose, however, that $G$ is a $T$--homogeneous
Lipschitz domain or, more generally, the boundary Harnack principle
holds on $T_0\cap G$: i. e., there is $\var>0$ and a constant $C>1$
such that for every pair of positive harmonic functions $u$ and $v$,
locally bounded and vanishing q. e.  near each point
of $\partial G\cap\{1-\var<|z|<1+\var\}$, we have
$$
\frac{u(\zeta)}{v(\zeta)}\leq C\frac{u(z_0)}{v(z_0)}
$$
for all $\zeta\in T_0\cap G$ and for some
$z_0\in T_0\cap G$. For instance, this happens if $T_0$
consists of finitely many arcs and in a neighborhood of each end-point
of these arcs the boundary of $G$ is a (possibly rotated) graph of a
Lipschitz function (cf. \cite {5, p. 178}).
By homogeneity, the same constant $C$  works on
each $T_n\cap G$. So, given a function $v\in \sP$, if
$H$ is the Martin function constructed in Theorem 0.5, let $z_0$ be
a point on $T_0\cap G$ where $M(1,H)$ (as defined in Proposition 4.1) is
attained. Then for all $z\in T_n\cap G$,
$$
v(z)\leq C\frac{H(z)}{H(T^nz_0)}v(T^nz_0)\leq
Cv(T^nz_0)\leq\frac{B}{\omega(z_0, T_n)}
$$
where we used (3.9) for the last inequality. Thus $\beta(v)<\infty$
(recall (3.2)) and so $v\in \sF$.

\subhead 6.\ Special properties of the spectrum
\endsubhead
Let $G$ be a $T$--homogeneous domain,
and $\r(G)$ be associated to $G$  as in Theorem 0.5.
(For example, as we noted after the statement
of Theorem 0.19, if $G=\{z; |\arg z| < \theta\}$, then
$\rho(G)={\pi}/{2\theta}$.) 
We first show that $\rho(G)$ is a (strictly) monotonic set function.

\proclaim{Theorem 6.1} Let $G_1\subset G_2$ be $T$--homogeneous
domains such that $E=G_2\backslash G_1$ has positive capacity.  Then
$\rho (G_1) > \rho (G_2)$.  \endproclaim

\demo{Proof} Since $TE=E$, it follows that 
$\cp(E\cap \{1\leq |z| < T\}) > 0$.
Without loss of generality, we may choose a compact set $K\subset
E\cap \{1 < |z| < T\}$ of positive capacity such that $K\subset\subset
G_2\cap \{1 < |z| < T\}$; if
$E\subset \{ |z|=1\}$, we replace $G_2,G_1$ by $\lm G_2,\lm G_1$ for $\lambda$
close to 1.

Let $T_n=\{ |z|=T^n\}$ and $z_0\in G_1\cap G_2$ with $|z_0|=1$.

Recall that for a $T$--homogeneous domain $G$ and a compact set $K\sbt
G$ the harmonic measure $\om (z,K,G)$ satisfies the equality  
$\om (Tz,TK,G)=\om (z,K,G).$ 

Now fix $n > 1$ and for $0\leq j\leq n-1$ put 
$$ m_j=\inf_{z\in T^j
K}\omega (z, T_n, G_2).\tag 6.2  $$ 
By the maximum principle, 
$$
\omega (z_0, T_n, G_2)\geq \omega (z_0, T_n,
G_2\backslash\bigcup^{n-1}_0 T^j K)+ \sum^{n-1}_0 m_j\omega (z_0,
T^j K,\ G_2\backslash \bigcup^{n-1}_{\ell=1} T^\ell K),
\tag6.3 $$
where $m_j$ is from (6.2) and since $K\subset\subset G_2$, the Harnack
inequality yields 
 $m_j\geq C\omega (T^j z_0, T_n, G_2)$ for some $C=C(G_2)>0.$
Hence we have from (3.15) (for a different $C$) that $$ m_j\geq
C\sup_{\zeta\in T_j}\omega (\zeta, T_n, G_2).\tag6.4 $$

Using the ideas of \S 3, let $I_j=\{ z\in T_j\cap G_{1}$,
$d(z,\partial G_1) > \var T^j\}$, where $0 < \var < d (z_0,\partial
G_1)$. Then
$$
\gathered
\omega (z_0, T^j K,\ G_2\backslash\bigcup^{n-1}_{\ell=0} T^\ell K)\geq
\omega (z_0, I_j, G_1)
\inf_{\zeta\in I_j} \omega(\zeta, T^j K,
G_2\backslash\bigcup^{n-1}_{\ell=1} T^\ell K)\\ 
\geq \omega (z_0, I_j, G_1)\inf_{\zeta\in I_j} \omega(\zeta, T^j K,
G_2\backslash \bigcup^\infty_{-\infty} T^\ell K)= A_j\omega (z_0, I_j,
G_1), 
\endgathered
$$
where, $A_j=A_0$ is independent of $j$.
Since $K\cap T_0=\varnothing$, $A_0 > 0$.

The argument which gave (3.10) shows that to  $\var < \var_0$
corresponds $n_0=n_0(z_0)$ so that 
$$
\omega (z_0, I_j, G_1)\geq {1\over 2}\omega (z_0, T_j, G_1)\quad(j>n_0).\tag6.5
$$
Hence (6.4) and (6.5) imply that when $j>n_0$,  each term in the sum
in (6.3) is greater than 
$$
{1\over 2} C A_0 \omega (z_0, T_j, G_1)\sup_{\zeta\in T_j} \omega
(\zeta, T_n, G_1) 
\geq C_1 \omega (z_0, T_n, G_1),
$$
if we take $C_1=C A_0/2$.
This transforms (6.3) to
$$
\omega (z_0, T_n, G_2)\geq (1+ (n-n_0) C_1)\omega (z_0, T_n, G_1),
$$
and so Theorem 6.1 follows from (4.4).        .
\qed
\enddemo

We can now prove  Theorem 0.16.
\demo{Proof of Theorem 0.16}  Let
$D\sbt \BT^2_P$ be connected on spirals and let $G$ be a component
of $\phi^{-1} (D)$. 
Proposition 0.10 implies that $G$
is a $T^{k}$-homogeneous domain, for some $k\in \Bbb N$.
Thus equation  (0.7) from Theorem 0.5 (with $T$ replaced by $T^{k}$) 
implies that $\r_0:=\r (G)$ is a positive point of
$Spec\ D$ with eigenfunction
$q(z)=H(e^z)e^{-\r(G)\Re z}.$
Together with Proposition 2.1 this proves
Theorem 0.16, (1).

Let $\r^*:=\r(D)$ be the minimal positive eigenvalue
of the boundary problem (0.15). It exists
because the set of positive eigenvalues is not
empty as we have just shown, is discrete
without any finite point of condensation
and does not contain zero
(Propositions 1.36, 1.37).
Now we are going to prove that $\r(G)=\r(D)$ and this
proves Theorem 0.16, (2).

Since $\r^*\leq\r_0$ we must prove $\r^*\geq\r_0.$

Denote by $q^*(z)$ the
eigenfunction  corresponding to $\r^*.$ We may
assume that $q^*(z_1)=1$ for some $z_1\in D$.

Let $G$ be a component of  $\phi^{-1} (D)$, and let 
$v^*(z)=q^{*} (\log z)|z|^{\rho},\ z\in G,$ so that $v^{*}$ and
$q^*$ are related by (0.8). Also, let  $G^*\sbt G$ be the component
of $\{v^*(z)>0\}$ which contains the preimage of $z_1.$ Then
$v^*$ is positive harmonic in $G^*$ and vanishes
quasi-every\-where on the boundary.
By Theorem 6.1 $\r^*\geq\r_0,$ and this establishes the remaining
assertion of Theorem 0.16. 
\qed
\edm
\demo {Proof of  Theorem 0.17} This follows
directly from Theorem 6.1 and the equality
$\r(G)=\r(D).$ Indeed,
$D_1\sbt D_2 \Longrightarrow G_1\sbt G_2$
for an arbitrary $G_2\in \phi^{-1}(D_2)$ and  $G_k,\ k=1,2$ corresponding to
$D_k,\ k=1,2$ by (0.9). The supposition
$\cp (D_2\setminus D_1)>0$ implies
$\cp (G_2\setminus G_1)>0$ since analytic maps
 preserve
positive capacity.
\qed
\edm
\proclaim {Proposition 6.6} Let
$D_n\uparrow D$ be a sequence of domains in
$\BT^2_P$ and $q_n,\ n=1,2,...$ are the
corresponding normalized by condition
$q_n(z_0)=1, z_0\in D$ solutions of the
problem (0.15).

Then $\r (D_n)\downarrow \r(D)$ and
$q_n\rightarrow q$ uniformly in any compact
set $K\sbt D$ and quasi-everywhere on $\p D.$
\ep
\demo {Proof}Since each $D_n$ can be approximated
from inside by smooth domains (see, e.g. \cite {12}), we can suppose
that each $D_n$ is smooth.

Set $\r^*:=\liml_{n\rightarrow \iy}\r(D_n).$
This limit exists because the sequence
$\r(D_n)$ decreases monotonically and is  bounded
below by $\r (D).$  Consider the sequence $\{H_n\}$ of
functions
$$H_n(z):=q_n(\log z)|z|^{\r (D_n)}.$$
Each $H_n$ is  positive harmonic  the  domain $G_n$
(which corresponds to $D_n$ by (0.9)), vanishes
on the boundary, and the sequence $\{H_n\}$ is compact.
Consider any convergent subsequence
$H_k\rightarrow H=q^*(\log z)|z|^{\r^*}.$

Since the $ \{q_n\}$ are normalized and converge
uniformly on compacta, we have
$q^*(z_0)=1.$
In addition, since $\{H_k\}$ (where each $H_k$ is extended
to be zero outside $G_k$) is a sequence
of subharmonic functions in $\BC,$ the function
$H$ is zero quasi-everywhere on $\p G$ by the theorem
of H. Cartan (\cite {10}, Chapter 7). It is also positive
harmonic in $G$.
By Theorems 0.5
and 0.16 $\r^*=\r (D)$ and $q^*=q.$ \qed
\edm

Let $G$ be a component of $\phi^{-1}(D)$ 
(see, (0.9)).  The point $0\in \p G$ plays a role analogous to that
of $\infty$, and this provides information which will supplement Proposition
1.37.
\proclaim{Proposition 6.7}
Let $D$ be a domain in $\BT^2_P$ which is connected on spirals so that
$Spec\ D$ is non-empty. Then the largest negative point in $Spec\ D$
is $-\rho(D)$.
\endproclaim
\demo{Proof}Suppose first that $D$ is smooth
enough and hence $\sF =\sP$ (see \S 5).
Let $G$ be related to $D$ as above
and let $H^+$ be the function provided by Theorem 0.5.
If $G(z,w)$ is the usual Green function for the domain $G$ and
$z_0\in G$ is a given base point, then
$$
H^+(z)=\lim_{n\rightarrow+\infty}\frac{G(z,T^nz_0)}{G(z_0,T^nz_0)}.\tag6.8
$$
Now recall that in Proposition 1.37 we defined the domain $D_-$ as
the $\bT^2_P$-domain obtained from $D$ via the map $z\mapsto -z$. This
corresponds to changing the plane-domain $G$ to $G_-$ via the map
$z\mapsto 1/z$. Define $H_{G_-}^+$ as above for the domain $G_-$. Then
letting $H^-(z)=H_{G_-}^+(1/z)$ we obtain from (6.8) that
$$
H^-(z)=\lim_{n\rightarrow+\infty}\frac{G(z,T^{-n}z_0)}{G(z_0,T^{-n}z_0)}.
$$
Let $\rho=\rho(D)>0$ and $\sigma=\rho(D_-)>0$. These constants
give
$$
H^+(Tz)=T^\rho H^+(z)\qquad H^-(T^{-1}z)=T^\sigma H^-(z).
$$
But by definition
$$
H^+(Tz)=\lim_{n\rightarrow+\infty}\frac{G(Tz,T^nz_0)}{G(z_0,T^nz_0)}\frac{G(z,T^nz_0)}{G(z,T^nz_0)},
$$
so that
$$
T^\rho=\lim_{n\rightarrow+\infty}\frac{G(Tz,T^nz_0)}{G(z,T^nz_0)}.
$$
Likewise we find that
$$
T^\sigma=\lim_{n\rightarrow+\infty}\frac{G(T^{-1}z,T^{-n}z_0)}{G(z,T^{-n}z_0)}.
$$
Notice that $z$ can be chosen arbitrarily. In particular if we let $z=z_0$:
$$
\aligned
T^\sigma  = &
\lim_{n\rightarrow+\infty}\frac{G(T^{-1}z_0,T^{-n}z_0)}{G(z_0,T^{-n}z_0)}
 =
\lim_{n\rightarrow+\infty}\frac{G(T^{n}z_0,Tz_0)}{G(T^nz_0,z_0)}\\
 = &
\lim_{n\rightarrow+\infty}\frac{G(Tz_0,T^nz_0)}{G(z_0,T^nz_0)}
 =
T^\rho:
\endaligned
$$
$\rho=\sigma$.

Let $D$ be an arbitrary domain connected on spirals.  There exists a
sequence $D_n$ of domains with smooth boundary such that $D_n\uparrow
D.$ Then the assertion of Proposition 6.7 for any domain $D$ follows by
Proposition 6.6 that is also a corollary of Theorem 0.5.  \qed \edm

Proposition 6.7 shows that $\rho(D_-)=\rho(D)$.  This
had been conjectured by V. Matsaev in a more general form,
namely, that
 $Spec(D)=Spec(D_-)$. In the generality of this paper, this question
seems to be open at the moment.

\subhead 7.\ Green function and Dirichlet problem\endsubhead
The following theorem defines the Green function corresponding to
$L_\rho$ and gives its properties. 
\proclaim{Theorem 7.1} Let $D\subset \Bbb T^2_P$ be a domain, 
$\r\notin Spec\ D$ be  given and
$L_\r$ as  in $(0.15)$. Then  there exists a fundamental
solution $g_{L_\rho} (z,\zeta,D)$ such that
$$L_\r g_{L_\r}(z,\z)=\delta_\z(z)\tag 7.2$$
$$\lim\limits_{z\rightarrow z'}
g_{L\rho} (z,\zeta, D)=0\tag 7.3$$
when $\zeta\in D$ and $z'\in \partial D$ is a
regular point in the sense of potential
theory, and the limit is uniform for $\z \in K$, $K$ compact in $D,$
and
$$ g_{L_\rho} (z,\zeta,D)\leq 0\qquad (z,\zeta\in D), \ 0< \r<\r(D).\tag 7.4$$
\endproclaim

\demo{Proof}From the Fredholm theorem we obtain 
that the equation 
$$\Cal G_D (\r):=I+2\r G^1_{D}+\r^2 G_{D}=f$$ has a unique solution
$q\in L^2(\Bbb T^2_P)$ for $\r\notin Spec$ and $f\in L^2(\Bbb T^2_P).$
Take $f:=g(z,\z,D),$ where $g$ is the Green function of the Laplace
operator and $g(\cdot,\z,D)\in L^2(\Bbb T^2_P)$ by Proposition 1.12
uniformly with respect to $\z$. Recall that $G^1_{D}q(z),$ $
G_{D}q(z)\ri 0$ as $z\ri z_0$ for every regular point of $\p D$ (see (1.31)).

 Thus $q$ satisfies (7.3). One can show, following the proof
of Proposition 1.28 , that (7.2) is satisfied too.

Let us prove (7.4). If $g_{L_\rho},$ which is a 
real function for a real $\r,$ were to change sign in $D$,
consider a component $$ D^- = \{z\in \Bbb T^2_P; g_{L_\rho}
(z,\zeta,D) < 0\}.  $$ Then $D\setminus D^-$ has positive capacity, and
$D^-\subset D$ is a domain in which $L_\rho g_{L_\rho}=0$.  According
to Theorem 0.17, $\rho (D^-) > \rho(D)$.  Since $\rho < \rho(D^-)$, we
have from the definition of $\rho (D^-)$ that $g_{L_\rho}\equiv 0$ in
$D^-$.  Note that $g_{L_\rho}$ can not be zero in $D$ without changing
sign in $D$ because of the maximum principle for 
the harmonic function $H(z)=g_{L_\rho}(z)
e^{\r z}.$
\edm
\demo {Remark}The explicit form of 
 $g_{L_\rho}$ is given by the expression
$$g_{L_\rho}(z,\zeta)=\sum_{k=-\iy}^{\iy} g (e^z, 
e^{\zeta+kP}, G)
e^{\rho(\xi+kP)} e^{-\rho x}\ (z=x+iy, \zeta=\xi+i\eta).
$$
where   $g(\cdot,\cdot,G)$ is the Green function of the initial domain $G\in \phi^{-1}(D)$ related to $D$
by (0.9). One can prove the convergence of 
this series for $0<\r<\r(D)$ and $z\neq \z,$ using (4.4).  
 \edm

Now consider the Dirichlet problem

$$L_\rho q=0,\ q|_{\p D}=f\tag 7.5$$ where $f$ is continuous in $\p
D,$ $D$ is regular domain (i.e., sufficiently smooth so that $q(z)\to
f(\z)$ while $z\to \z \ \forall \z\in \p D$), $\r\notin Spec\ D.$

\proclaim{Proposition 7.6} If $0< \r\notin Spec\ D$, then the solution
$q(z,f,D)$ to the problem $(7.5)$ is unique.  If $0<\r<\r(D)$ and $f
\geq 0$ on $\partial D$, then $q \geq 0$ in $D$.  \endproclaim
\demo{Proof} If $q_1$ and $q_2$ solve (7.5), then $q_1-q_2$ would be
the unique solution to the homogeneous problem (0.15), and so by our
assumption on $\rho$, $q_1\equiv q_2$.

Next, let $q$ solve (7.5) where $f\geq 0$ on $\partial D$, and suppose
$q(z_0) < 0$.  Then $q=0$ on the boundary of a connected component of
the open set $D^-=\{z;\, q(z) < 0\}$.  But $D^-\subset D$, so
$\rho(D^-) \geq \rho (D) > \rho$.  Thus $q\equiv 0$ in $D^-$, which is
a contradiction.  \qed \edm Let $f$ be upper semicontinuous on
$\partial D$, and consider a sequence of continuous functions
$f_n\downarrow f$.  The corresponding sequence $q_n=q(z,f_n,D)$ (using
$f_n$ in (7.5)) converges monotonically, and defines a unique solution
to (7.5) for upper semicontinuous $f$. The same holds for lower
semicontinuous functions. Since every measurable function can be
represented as a sum of functions of these two types, the solution is
defined and unique for all mesurable functions. It can be equal to
$\iy$ or $-\iy$.  This is a generalized solution in the sense of
Wiener to the problem (7.5); see, for example, \cite {10}.

\subhead 8.\ Subfunctions with respect to the operator $L_\rho$\endsubhead

An upper semicontinuous function on $\Bbb T^2_P$ is called an $L_\rho$--{\it
subfunction} if, $L_\rho v\geq 0$ in the sense of distributions.  If
both $v$ and $-v$ are $L_\rho$--s.f., we call $v$ an
$L_\rho$--{\it function}.  The theory of $L_\r$-subfunctions parallels
that of subharmonic functions because of \proclaim{Proposition 8.1} A
function $v$ is an $L_\rho$--s.f. iff the function
$$
V(z)=v(\log z) |z|^\rho
$$
is subharmonic in $\bC$ and satisfy
$$
V (e^P z) e^{-\rho P}=V(z).\tag8.2
$$
\endproclaim
\demo{Proof} That $V$ is well--defined,upper--semicontinuous  and
satisfying (8.2)  follows from properties $V$ inherits from $v$.

We claim that $\Delta V\geq 0$ in $\Cal D' (\bC)$.
If $\Psi\in\Cal D (\Bbb C \backslash 0)$, then $\Psi$ may be written
as $\Psi(z)=|z|^{-\rho}\psi (z)$, with $\Psi\in\Cal D(\Bbb C\backslash
0)$, and we may suppose that the support of $\Psi$ is contained in a
sector 
$\Delta(\alpha,\beta,R,P)=\{re^{i\varphi};\varphi\in(\alpha,\beta),
|\log (R/r)| < P\}$. 
Let
$$
\align
& r=e^x,\ \varphi=y,\ z=x+iy,\ \hat\psi (z)=\psi (e^z), v(z)=e^{-\rho
x} V(e^{z}),\\ 
& \Delta (\psi (z) |z|^{-\rho})=L_\rho^*\hat\psi (\log z) e^{-(\rho-2)x},
\endalign
$$
where $L_\rho^*=L_{-\r}$ is symmetric 
to $L_\rho$ .

Then
$$
\gathered
\langle\Delta\psi (\cdot) |\cdot|^{-\rho}, V\rangle =\int
L_\rho^*\hat\psi (z) v(z) dx dy\\ 
=\langle\hat\psi, L_\rho v\rangle_{\Bbb T^2_P}\geq 0,
\endgathered
$$
since $L_\rho v\geq 0$ in $D'(\Bbb T^2_P)$.
Thus $V$ is subharmonic in $\Bbb C\backslash 0$, and since $V$ is
bounded near 0, $V$ extends to be subharmonic at the origin. 

The sufficiency follows in the same way.
\qed
\enddemo

\proclaim{Corollary 8.3}The following holds
\itemitem{(1)}If $v_1,\dots, v_k$ are $ L_\rho-s.f$ then
$\max\limits_i v_i$ is an $L_\rho$--s.f.; 
\itemitem{(2)}The set of $L_\rho$ subfunctions form a positive cone;
\itemitem{(3)}The set of $L_\rho$ subfunctions are closed under
translation of coordinates; 
\itemitem{(4)}If $d\mu\geq 0$ and $v$ is $L_\rho$--s.f., then
$$
\int_{\Bbb T_P^2} v(z-\zeta)d\mu(\zeta)
$$
is an $L_\rho$--s.f.
\endproclaim

Let $v$ be an $L_\rho$--s.f. Consider the $C^\infty$-- function in
$\Bbb C\setminus 0$ 
$$
V_\var (z,v)=
|z|^{-2}\int\alpha_\var (\zeta/z ) v(\log\zeta)|\z|^{\rho}
d\z
$$
where $\alpha_\var\in\Cal D(\Bbb C),\ \alpha_\var\geq 0,
\ \alpha_\var (z)=0\ $ for $|z-1|>\var $ and
$\int\limits_{|\zeta -1| < \var}\alpha_\var(\zeta)
 d\zeta=1. $
It is easy to verify that $V_\var$ is subharmonic and
satisfies (8.2). Thus
$$
\Cal M_\var(z,\cdot):=
e^{-\rho x}V_\var (e^z,\cdot)
\tag 8.4$$
is an $L_\rho$--s.f., and a straightforward computation shows that
$$\lim\limits_{\var\to 0}
\Cal M_\var (z, v)=v(z).\tag 8.5$$

This will yield the first assertion of 
\proclaim{Proposition 8.6} $(1)$ Every $L_\rho$--s.f.~is a decreasing
limit of a sequence of 
infinitely differentiable $L_\rho$--s.f.'s.

$(2)$ A non-zero $L_\rho$--s.f. $v$ cannot attain a local non-positive
maximum in a domain  $D\subset \Bbb T^2_P$.  \endproclaim
\demo{Proof}

Now we prove (2).
Let $v(z_{\max}):=-c\ (\leq 0)$ be the maximal value of $v(z).$ We
apply the maximum  
principle to the subharmonic function $V$ associated to
$v_1(z):=v(z)+c$ in Proposition 8.1 and  
obtain that $V(z)\leq 0$ and $V(z_{\max})=0,$ 
hence $V(z)\equiv 0.$ Thus $v(z)\equiv -c.$ If 
$c>0,$ then $L_\r(-c)<0$ that is a contradiction. 
\qed \enddemo

\subhead 8.7. Green Potential\endsubhead Let $D\subset \Bbb T^2_P,\ 0<\rho <
\rho(D)$, $v$ an $L_\rho$--s.f., and let $g_{L_\rho}$ be the Green
function of $L_\rho$, cf. \S7.  Consider the {\it Green potential} of a
measure $\nu$ on $D$:
$$ \Pi (z,\nu)=\int_D g_{L_\rho} (z,\zeta,D)
d\nu(\zeta).  $$

\proclaim{Proposition 8.8} $\Pi(z,\nu)$ is an $L_\rho$--s.f.~in $D$ with
$$
L_\rho \Pi (\cdot,\nu)=\nu\qquad\text{ in }\qquad\Cal D'(\Bbb T^2_P).
$$
Moreover,  if $supp(\nu) \sbt\sbt D$
or if $\nu$ has bounded density in a neigborhood of a regular boundary point
$z_0\in\p D,$
then
$$ \lim\limits_{D\ni z\rightarrow z_0}\Pi (z,\nu)=
0.$$
\endproclaim

\demo{Proof} The function 
$\Pi(z,\nu)$, being the potential of a negative kernel, is upper
semicontinuous.
Next, let $\Psi > 0\in\Cal D(\Bbb T^2_P)$.
By Theorem 7.1, we have
$$\align
\langle L_\rho \Pi,\psi\rangle &=
\int L^*\psi\int_D g_{L_\rho} (z,\zeta,D)d\nu(\zeta) 
=\int_D d\nu (\zeta)\int L_z^*\psi g_{L_\rho} (z,\zeta) dz\\
&=\int_D \psi(\zeta)d\nu(\zeta) > 0,
\endalign$$
so that $L_\rho\Pi(\cdot,\nu)=\nu$ in $\Cal D'(\bT^2_P).$

Let  $\nu'$ be  the density of $\nu.$ It is straightforward that for
any $z_0\in D\cup \p D$
$$\sup_{|z-\z_0|<\eps}\int_{|\z-z_0|<\eps }
|g_{L_\rho} (z,\z,D)|d\z
\leq C\int_{|\z-z_0|<\eps }
|\log |\z-z_0||d\z$$
with a constant $C$ which does not depend on $\eps$ and $z_0.$
Hence,
$g_{L_\rho} (z_n,\z,D)\nu'(\z)$ is a Lebesgue sequence when $z_n\ri z_0$.
\qed\enddemo

Let $v$ be an $L_\rho$--s.f.  Then since $v$ is upper semicontinuous,
the solution $q(z)$ of the problem (7.5) with boundary data $v$ is
defined for any regular domain $D$ for which $0<\rho < \rho (D)$; this
function $q=q(z)=q(z,f,D)$ is called {\sl the least} $L_\rho$--{\sl
majorant} of $v$ in $D$.

\proclaim{Proposition 8.9 (Riesz Theorem)} Let $v$ and $q$ be as
above, and $0<\rho < \rho (D)$.  Then
$$
v (z)=q(z) + \Pi (z,v),
$$
where $\nu=L_\rho v$ and $q$ is the least $L_\rho$--majorant of $v$ in $D$.
\endproclaim

\demo{Proof} If $v$ is sufficiently smooth, then $q=v-\Pi (z,\nu)$ is
an $L_\rho$--function which agrees with $v$ on $\partial D$
quasi--everywhere.  For an arbitrary $v$, we obtain this
representation for a sequence $v_n\downarrow v$ of smooth
$L_\rho$--s.f.'s, and take the limit in $\Cal D'(\Bbb T^2_P)$.\qed

\enddemo
It is clear that $q(z,v,D)$ is a majorant of
$v.$ The next assertion, which follows from Proposition 7.6,
 shows that it is
the least majorant.

\proclaim{Proposition 8.10} Let $ D\subset \Bbb T^2_P$ be a domain and
$0<\rho < \rho (D)$.  Let $f$ be upper semicontinuous on $\partial D$,
and let $v$ be an $L_\rho$--s.f.~in $\Bbb T^2_P$.

Then if $f\geq v$ on $\partial D$, it follows that $q=q(z,f,D)\geq v$
in $D$.  \endproclaim
\proclaim{Theorem 8.11 (Sweeping)} Let $\rho < \rho(D)$, $v$ an
$L_\rho$--s.f.~and $q(z,v,D)$ the least $L_\rho$--majorant of $v$.
Then the function
$$
v(z,v,D):=\cases v(z)&\text{$(z\in \Bbb T^2_P\backslash D)$}\\
q(z,v,D)&\text{$(z\in D)$}.\endcases
$$
is an $L_\rho$--s.f.
\endproclaim

\demo{Proof}
We need check this only in a
neighborhood of $\p D$ where it follows from
the inequality
$$
v(z)\leq q(z,v,D)\qquad (z\in D)
$$
and Corollary 8.3(1).
\qed
\enddemo

\example{Example 8.12} Consider the situation (0.2).We consider $v$ in
the sector $G=\{z=re^{i\theta}: \alpha < \theta <\beta\}$ and
associate to $G$ the subdomain $D=\{z\in \BT^2_P: \alpha \leq
y\leq\beta\}$.  Since $\rho (G)=\pi/(\beta-\alpha)$, we have that
$\rho (D)=\pi/(\beta-\alpha)$ independent of $P$, and the condition
$\rho <\rho (D)$ reduces to the classical requirement $(\beta-\alpha)
< \pi/\rho$.  If $h$ satisfies (0.1), then the smallest
$\rho$--trigonometric majorant of $h$ on $(\alpha,\beta)$ is
$$
H(\varphi)={h(\alpha)\sin\rho(\beta-\varphi)+h(\beta)\sin\rho(\varphi-\alpha)
\over \sin\rho (\beta-\alpha) },
$$
and Theorem 8.11 reduces to the so--called fundamental relation for
$\r$--t.c.functions \cite {16, (1.69)}:\ if $\max\limits_{i,j}
|\varphi_i-\varphi_j| < \pi/\rho,\ i,j=1,2,3$, then
$$
h(\varphi_1)\sin\rho(\varphi_2-\varphi_3)+h(\varphi_2)\sin\rho(\varphi_3-\varphi_1)+h(\varphi_3)\sin\rho(\varphi_1-\varphi_2)\leq 
0.\tag 8.13
$$
\endexample

The procedure of constructing the least harmonic majorant $H(z,u,G)$
associated to a subharmonic function $u$ in $G\subset\Bbb C$ is called
the `sweeping of the masses' of $u$.  Thus we may call our
construction of $v(\cdot,v,D)$ the `sweeping' of $L_\rho$--masses of
$v$.

\subhead{8.14 Representation Theorems}\endsubhead It is natural to
describe new classes of functions in terms of independent parameters.
For example, if $\rho$ is not an integer, a $\rho$--trigonometrically
convex function $h$ has the representation
$$
h(\varphi)={1\over
2\rho\sin\pi\rho}\int_0^{2\pi}\cos^*\rho(\varphi-\psi-\pi)d\Delta(\psi)
$$
where the function $\cos^*\rho\varphi$ is the $2\pi$--periodic
extension of $\cos\rho\varphi$ from $(-\pi,\pi)$, and
$\Delta(d\psi)=[h''(\psi)+\rho^2h(\psi)]d\psi$ is a positive measure
(see \cite {16, Ch.~1, Theorem 24} and \cite {6, Ch.~1}, where the case
$\rho\in\Bbb N$ is also considered).  From our point of view,
$d\Delta$ is the independent parameter for the class of $\rho$--
t.c. functions, and its connection with the zero--distribution of
functions of completely regular growth (Levin-Pfluger functions) is
the theme of \cite {16, Ch.~2}.

We now consider the situation corresponding to the operator $L_\rho$
and on $\Bbb T^2_P$. 
Propositions 8.15, 8.17 and 8.19 generalize
\cite {16, Theorem 24}. The proof is new also
for $\r$--t.c.functions.

First, let $\rho\notin \Bbb Z$ and $E_\rho$ the fundamental solution
of $L_\rho$ on $\Bbb T^2_P$, as in Proposition 1.1.  For a measure
$\nu$, consider the potential $$ \Pi_\rho (z,\nu)=\int_{\Bbb T^2_P}
E_\rho (z -\zeta) d\nu(\zeta).  $$ As in the proof of Proposition 8.8,
$\Pi_\rho$ is an $L_\rho$--s.f.~and $L_\rho \Pi_\rho=\nu$.
\proclaim{Proposition 8.15} Let $\rho > 0,\ \rho\not\in \Bbb Z$.  Then
every $L_\rho$--s.f.~$v$ on $\Bbb T^2_P$ may be represented as $$
v(z)=\Pi_\rho (z,\nu), $$ where $\nu=L_\rho v$.  \endproclaim

\proclaim{Lemma 8.16} Let $L_\rho q=0$ on $\Bbb T^2_P$.
Then $q\equiv 0$ when $\rho\not\in \bZ,$ and $q(z)=\Re\{C e^{i\rho
y}\},\ C\in \Bbb C,$ when  $\rho\in\bZ$. 
\endproclaim

\demo{Proof} Exactly as in Proposition 1.1, we see that the Fourier
coefficients 
$\{q_{m,k}\}$ of $q$ must be chosen so that $$ q_{m,k}\left[-\left({2\pi m\over
P}\right)^2-k^2-2\rho i {2\pi m\over P}+\rho^2\right]=0, \
(m,k)\in\bZ^2.  $$ When 
$\rho\not\in\bZ$, this forces all $q_{m,k}$ to vanish, and when
$\rho=p\in \BZ$, the bracketed term vanishes when $m=0,\ k=\pm p.$
  In this
case, if $q_{0,\pm p}\not= 0$, we have $q=c_1 e^{-ip y} + c_2
e^{ipy}$, and since $q$ is real, the Lemma follows.
\qed
\enddemo

\demo{Proof of Proposition 8.15} We apply $L_\rho$ to $q\equiv
v-\Pi_\rho$, and note that $q$ 
is an $L_\rho$--function on $\Bbb T^2_P$. Hence
$q\equiv 0$ as follows from Lemma 8.16.
\qed
\enddemo

Consider now the case $\rho=p.$ The following
assertion generalizes the corresponding condition for the
Phragm\'en-Lindel\"of indicator \cite {16, 1.81}. It  describes a
symmetry of mass distribution in the case of 
integral $\rho.$
\proclaim{Proposition 8.17} Let $\rho=p\in
\Bbb Z,\ p\geq 1,$
$v$ be an $L_p$--s.f. on $\Bbb T^2_P$ and
$\nu = L_p v$. Then
$$\int_{\Bbb T^2_P} e^{\pm i p y} d\nu (z)=0.
\tag 8.18$$
\endproclaim
\demo {Proof} The functions $ e^{\pm i p y}\in
\Cal D(\Bbb T^2_P)$ and are solutions to the
equation $  L^*_p q=0$ on $\Bbb T^2_P.$ Thus
$$\int_{\Bbb T_P} e^{\pm i p y} d\nu (z)=\langle e^{\pm i p \cdot},
L_p v\rangle= 
\langle L_p^* e^{\pm i p\cdot},v\rangle=0.\qed$$
\enddemo

Let $E_p'(z)$ be arbitrary generalized fundamental
solution from Proposition 1.10. Set
$$\Pi'_p (z,\nu):=\int_{\Bbb T^2_P}
E_p' (z-\zeta)d\nu(\zeta).$$
The potential is defined uniquely because of
Propositions 1.10 and 8.17.

The next assertion generalizes the
representation of the Phragm\'en-Lindel\"of
indicator for functions of  integral order
\cite {16, 1.82}

\proclaim {Proposition 8.19} Let $p$ be an integer and
$v$ be an $L_p$--s.f.~on $\Bbb T^2_P$. Then
$$v(z)=\Pi_p'(z,\nu) +\Re(Ce^{ipy}),\tag8.20
$$
where $\nu=L_p v$, and $C$ is a complex scalar.
\endproclaim
\demo {Proof}Using Proposition 1.10 we have
$$L_p\Pi'_p (z,\nu)=\nu-\int_{\Bbb T^2_P}\cos
p\Im(z-\zeta)\,d\nu(\zeta).\tag 8.21$$ 
Hence Proposition 8.17 gives
$L_p\Pi'_p (z,\nu)=\nu.$ Thus the function
$q:=v-\Pi'_p (z,\nu)$ satisfies $L_pq=0$ on
$\Bbb T^2_P.$ Lemma 8.16 gives that
$q=\Re(Ce^{ipy})$
\qed
\enddemo

\subhead 9.\ $L_\rho$--subminorants\endsubhead An $L_\rho$--s.f.~$v$
is called an $L_\rho$--subminorant $(L_\rho$--s.m.) of a real--valued
function $m(z)$, $z\in D\subset \Bbb T^2_P$, if $$ v(z)\leq m(z),\ z\in
D.\tag9.1 $$ An $L_\rho$--subminorant $v_0$ is called the {\it
maximal} $L_\rho$--subminorant of $m$, if the conditions $\{w$ is an
$L_\rho$--subminorant$\}$ and $\{w\geq v_0\}$ imply that $w=v_0$ in
$D$.

\proclaim{Theorem 9.2} If $m(z)$ is continuous and has an
$L_\rho$--subminorant, then it has a unique maximal
$L_\rho$--subminorant. 
\endproclaim

\demo{Proof} This assertion follows by word--word repetition of the
proof for the case of maximal subharmonic minorant \cite {14}, which
we briefly sketch. 
The set of the subminorants is a partly
ordered set, because the semicontinuous regularization of the supremum
of any set of subminorants is also a subminorant.  Hence, there exists
a unique maximal element.
\qed
 \enddemo

For some applications it is desirable to have maximal
$L_\rho$--s.m.~when $m(z)$ need not be continuous.  At present, this
is possible only in certain cases, even in the classical case of the
Laplace operator.  First, let $m(z)$ be upper semicontinuous, and let
the sequence of continuous functions $m_n\downarrow m(z)$.  The
corresponding sequence of maximal $L_\rho$--s.m.'s decreases
monotonically and thus converges to an $L_\rho$--s.f.~$v_0$ which is
the maximal $L_\rho$--s.m.~for $m(z)$ as can readily be verified.

When $m(z)$ is not upper semicontinuous, we can, of course, consider
its upper semicontinuous regularization
$$
m^*(z)\colon =\lim_{\epsilon\to 0} \sup \{ m(z)\colon |z-\zeta| <
\epsilon\}\tag9.3 
$$
which is an upper semicontinuous function, and thus construct an
$L_\rho$--s.m.~for $m^*$.  However, we cannot ensure that the maximal
$L_\rho$--s.m.~of $m^*$ does not exceed $m$ for all $z$.

We present some positive results.

\proclaim{Theorem 9.4} Let $m=m_1-m_2$ where $m_1,m_2$ are $L_\rho$--s.f.'s,
and assume that $m$ has an $L_\rho$-subminorant.
Then there exists a unique maximal $L_\rho$--s.m.~$v^*_{\sup}$.
\endproclaim

\def\sM{\Cal M} \demo{Proof} Since $m$ has an $L_\rho$--s.m.~$v$, we
recall the notation (8.5) and observe that $$ \sM_\epsilon (z,m)\equiv
\sM_\epsilon (z,m_1)-\sM_\epsilon (z,m_2) $$ has the $L_\rho$--s.m.
$\sM_\epsilon (z,v)$ and so, by Theorem 9.2, has the unique maximal
$L_\rho$--s.m.~$v(z,\epsilon)$.  Then $$ v_{\sup} (z)\equiv
\limsup_{\epsilon\to 0} v(z,\epsilon)\leq m(z).  $$ In this
inequality, we refer to (9.3), and note that $v^*_{\sup} (z)$ is an
$L_\rho$--s.f.~that coincides with $v_{\sup} (z)$ everywhere except
perhaps on a set of zero capacity. This follows by Cartan's theorem
(\cite {10}, 
Ch.7) applied to the sequence of subharmonic functions
$u(z,\eps):=v(\log |z|,\eps)|z|^\r.$ In general, $v^*_{\sup}$ can
exceed $v_{\sup} (z)$.  However, under our special  hypotheses here, we claim
that $$ v^*_{\sup} (z)\leq m(z)\tag9.5 $$ everywhere.  Indeed, $m_2(z)
+ v^*_{\sup} (z)\leq m_1(z)$ outside of a set of zero capacity, so
$$ \sM_\epsilon (z,m_2)+\sM_\epsilon (z,v^*_{\sup})\leq \sM_\epsilon
(z,m_1).  $$ The formula (8.5) now gives (9.5).

We show that $v^*_{\sup}$ is the maximal $L_\rho$--s.m.
If not, there would exist an $L_\rho$--s.m.~$v_1$ exceeds $v^*_{\sup}$ on a set of positive measure (otherwise they coincide); thus we would have for some $z$ and $\epsilon$
$$
v^*_{\sup} (z) < \sM_\epsilon (z,v_1)\leq v(z,\epsilon),
$$
and this contradicts the definition of $v^*_{\sup}$.
\qed
\enddemo

\proclaim{Proposition 9.6} Let $v$ be the maximal $L_\rho$--s.m.~of
$m(z)$ such that 
$
v (z) < m(z)
$
on an open set $U$.
Then $v$ is an $L_\rho$--function in $U$: $L_\rho v=0.$
\endproclaim

This proof parallels that for a subharmonic function (see \cite {14}),
but we need some technical details.

\demo{Proof of Proposition 9.6}
Note that in a  disc $D_\dl:=\{z\in \BT^2_P:|z-z_0|<\dl\}$ such that
$D_\dl\cap TD_\dl=\emptyset$ we have 
$$q(z,v,D_\dl)=e^{-\r x}\int\limits_{|\z-z_0|=\dl} P(z,\z,\dl)
e^{\r \xi}v(\z)|d\z|$$
where $P(\cdot,\cdot,\cdot)$ is the Poisson kernel, $|d\z|$ is the 
element of length.
Since  
$$\psi(\dl):=\max\limits_{|z-\z|\leq \dl}
|e^{-\r( z-\z)}|-1=o(\dl),$$
we obtain that 
$$q(z,v,D_\dl)\leq (1+o(\dl))\max
\limits_{|\z-z_0|\leq \dl} v(\z).\tag 9.7$$

Suppose that $v(z_0)<m(z_0)$ and $v$ is not an $L_\r$ -function in a
neighborhood $U\ni z_0.$ Set $m(z_0)-v(z_0):=d.$ Choose a disc $D_\dl$
such that $\nu_v(D_\dl)>0$ and so small that $$(1+o(\dl))\max_{\z\in
D_\dl}v(\z)<v(z_0)+d/2,\tag 9.8$$ where $o(\dl)$ is from (9.7). The
inequality (9.8) is possible because of upper semicontinuity $v(z).$

 Let us replace $v$ in $D$ by its least $L_\r$ -majorant,
i.e. construct the sweeping $v(z,v,D)$ from Theorem 8.11.  Then
$v(z,v,D)<m(z)$ for $z\in D$ and hence for all $z\in \Bbb T_P^2.$ But
$v(z,v,D)> v(z)$ in $D_\dl$ because of the Riesz Theorem 8.9. Thus $v(z)$
is not the maximal minorant. This is a contradiction.  \qed \edm

\proclaim{Proposition 9.9} If $m(z),
\ z\in\Bbb T_P^2,$ is continuous,
the maximal $L_\rho$--s.m.~is continuous.  \endproclaim

Note from Proposition 9.6 that there is no problem at points where the
maximal $L_\rho$--s.m.~does not strictly exceed $m$. 

We are going to use the following

\proclaim{Theorem 9.10} Let $m(z), z\in\bC$ be continuous and have a
subharmonic minorant in $\bC$.  Then its maximal subharmonic minorant
is continuous.  \endproclaim

This fact was not obvious for us and we could not find a proof.
Thus we thank Prof.~A.~Eremenko for the following argument:

\demo{Proof} We prove continuity at $z=1$.  Let $m$ be the continuous
function and $u$ its maximal subharmonic minorant.  Since $u$ is
already upper semicontinuous, we need only show that for every
$\epsilon > 0$
$$
u(z) > u (1)-\epsilon\tag 9.11
$$
in some neighborhood of $z=1$.  Let $v$ be the sweeping of $u$ in a
neighborhood $U$ of 1 (in a small disc).  Then it is easy to see that
$u(z)\leq v(z) < m(1) + \epsilon/4 < m(z) + \epsilon/2$ in $U$.  Hence
$v-\epsilon$ is a subharmonic minorant of $m$, and so $u > v-\epsilon$
everywhere.  Since $v$ is continuous in the disk, we can find a
neighborhood of $z$ in which $v(z) > v(1)-\epsilon/2$.  Thus (9.11)
holds in this neighborhood.  \enddemo

\proclaim{Proposition 9.12} Let $m$ be continuous in $\BC$ and satisfy
the condition $(8.2)$.  Then its maximal subharmonic minorant also
satisfies $(8.2)$.  \endproclaim

\demo{Proof} Let $u(z)$ be a subharmonic minorant of $m(z)$.
Then

$$
u_1(z)\colon = [\sup_{n\in\bZ} u(e^{nP} z) e^{-\rho n P}]^*,
$$
where $[\cdot]^*$ is defined by (9.3), is a subharmonic minorant of
$m(z)$, $u_1(z)\geq u(z),\ z\in \BC$, and $u_1$ satisfies (8.2). 

\qed
\enddemo

\demo{Proof of Proposition 9.9} Set $m_1(\lambda)\colon
=m(\log\lambda)|\lambda|^\rho$.  It is continuous in $\bC$ and
satisfies the assumption of Proposition 9.12.  Thus its maximal
subharmonic minorant $v_1(\lambda)$ is continuous and satisfies (8.2),
so by Proposition 8.1, $v(z)\colon = v_1 (e^z) e^{-\rho x}$ is the
$L_\rho$--s.m.~of $m$. In particular $v$ is continuous.
\qed
\enddemo

If $m(z)$ is not continuous or even upper semicontinuous, its maximal
$L_\rho$--s.m.~may still be continuous.

\proclaim{Theorem 9.13} Let $m=m_1-m_2$ where $m_1,m_2$ are $L_\rho$--s.f.'s.
Then the maximal $L_\rho$--s.m.~is continuous if $m_1$ is continuous.
\endproclaim

\demo{Proof} Let $v$ be the maximal $L_\rho$--s.m.~of $m$, and let
$\nu_m=\nu_{m_1}-\nu_{m_2}$, where $\nu_{m_i}=L_\rho m_i,\ i=1,2$.
Set $g(z)=e^{\rho x}[m_1-m_2-v](z)$ and let $E=\{z:g(z)=0\}$. On
$\Bbb C\setminus E$ we have $g(z)>0.$

Now we use the following assertion (Grishin Lemma) \proclaim {Theorem
AFG \cite {7}} Let $g$ be a nonnegative $\dl$ -subharmonic function
and $\nu_g$ be its charge. Then the restriction $\nu_g|_E$ to the set
$E:=\{z:g(z)=0\}$ is a measure.  \endproclaim Thus
$$
\nu_v\leq \nu_{m_1}-\nu_{m_2}\leq \nu_{m_1}.
$$
on E. By Proposition 9.6  $\nu_v=0$ outside $E.$
Hence this also holds in $\Bbb C.$

It follows from Jensen's theorem that a subharmonic function $u$ is
continuous at a point $z_0$ if and only if
$$\int\limits_0^\epsilon
\frac {\mu_{u}\{z:|z-z_0|<t\}}{t}dt=o(1),\eps\ri 0.$$
By hypothesis $m_1$ and, hence, $u:=e^{\r x}m_1$
are continuous. Thus $e^{\r x}v$ is continuous,
and hence so is $v$.
\qed
\edm

\subhead 9.14. A new set characteristic\endsubhead Let $D\subset\Bbb
T^2_P$ be a domain, $\rho(D)$ be its order, and set $$ \lambda
(D)={1\over \rho (D)}.  $$ This characteristic is a ``natural''
monotonic functional and is zero on any domain which is not connected
on spirals.  We extend $\lambda$ to arbitrary sets in the standard
way.  If $D\subset \Bbb T^2_P$ is open, define $\lambda
(D)=\max\limits_i \lambda (D_i)$, where $\{D_i\}$ are the connected
components of $D$.  For a closed set $K\subset\Bbb T^2_P$ we define
$\lambda$ as $$ \lambda(K)=\inf_{D\supset K}\lambda(D), $$ where $D$
is open.  Finally, for any set $E\subset\Bbb T^2_P$ let $$
\overline\lambda (E)=\inf_{D\supset E}\lambda (D);\ \underline\lambda
(E)=\sup_{K\subset E} \lambda(K), $$ where $\{D\}$ are open and
$\{K\}$ are closed.

It is  important to know if a given function $m(z)$
has an $L_\rho$--s.m.~and  then describe its maximal
$L_\rho$--s.m.~(cf.\cite {3, 4}):

\proclaim{Theorem 9.15} Let $m(z)$ be a function on $\Bbb T^2_P$, and let
$$
E^+ (m)=\{z\in \Bbb T^2_P, m(z) > 0\}.
$$
If $m$ has a non-zero $L_\rho$--subminorant, then
$$
\overline\lambda (E^+ (m))\geq {1\over\rho}.
$$
\endproclaim

\proclaim{Theorem 9.16} Let $m(z)\geq 0$ be a continuous function and
$$
\lambda (E^+ (m)) > {1\over\rho}.
$$
Then $m(z)$ has a non-zero $L_\rho$--subminorant.
\endproclaim

 Theorem 9.15 follows directly from
\proclaim{Theorem 9.17} Let $v(z),\ z\in
\Bbb T^2_P$ be an $L_\rho$--s.f.~and $E^+ (v)$ be defined as in the
statement of Theorem 9.15. 
Then $\bar{\lambda}(E^+ (v))\geq 1/\rho$ or
$v\equiv 0$ in $\Bbb T^2_P.$
\endproclaim

\demo{Proof} Suppose the theorem false, and choose an open set
$D\supset E^+(v)$ such that $\lambda(D) <1/ \rho.$ Hence for each
component $D_i$ of $D,$ $\rho(D_i) > \rho$.  Since the function
$q\equiv 0$ is the unique solution of the boundary problem (0.15) in
each $D_i$ we obtain from Theorem 8.10 that $v(z) \leq 0,\ z\in
D_i$. Thus $v(z)= 0, \z\in \Bbb T^2_P$ by Proposition 8.6(2). \qed
\enddemo

For the proof of Theorem 9.16 we need
\proclaim {Lemma 9.18} For every domain $D$ with
$\rho (D)<\rho$ there exists a domain
$D(\rho)\sbt\sbt D$
 such that $\rho(D(\rho))=\rho.$ \endproclaim This follows from
Theorem 0.17 and Proposition 6.6.  \demo{Proof of Theorem 9.16} Since
$E^+(m)$ is open, the condition $\lambda(E^+(m)) > \rho^{-1}$ implies
there is a connected component $D\subset E^+ (m)$ with $\rho(D) <
\rho$.  Let $D_1:=D(\rho)$ be from Lemma 9.18 and $v_1$ be a solution
to (0.15).  Set $$ v(z)=\cases Cv_1(z),&\text{$z\in D_1$}\\
0,&\text{$z\not\in \Bbb T^2_P\backslash D_1$}, \endcases $$ with
$C<\min\limits_{z\in D_1}m(z)/v_1(z).$ Then $C v_1 (z) < m(z)$ for
$z\in D_1$.  We thus obtain an $L_\rho$--s.m.~of $m(z)$.  \qed

\enddemo

Here is another necessary condition.

\proclaim{Proposition 9.19} Let $m(z),\ z\in\Bbb T^2_P$ have an
$L_\rho$--s.m. ($\rho>0$). Then $$ \int_0^{2\pi} m(x+iy) dy\geq 0,\
\forall x.\tag9.20 $$ \endproclaim

\demo{Proof} Let $v(z)$ be an $L_\rho$--s.m., and associate to $v$ the
subharmonic function $V(z)=v(\log z) |z|^\rho$.  Since $V(0)=0$ we
have that $\ds\int\limits_0^{2\pi} V (r e^{i\phi}) d\phi\geq v(0)=0,\
\forall r > 0$.  Therefore (9.20) holds for $v$ and hence for $m$.
\enddemo

\subhead 9.21. Minimality \endsubhead An $L_\rho$--s.f.~$v$ is called
{\it minimal} if $v-\var$ does not have an $L_\rho$--s.m.~for any
$\var > 0$, (see also \cite {3, 4}).  We shall have many examples of minimal
$L_\rho$--s.f.~in $\Bbb T^2_P$ once we establish the following
theorem.

Denote by $\sH_\rho(v)$ the maximal open set in which $L_\rho v=0$;
i.e., $v$ is an $L_\rho$--function in $\sH_\rho(v)$.

\proclaim{Theorem 9.22=0.20} If there exists a connected component
$D\subset \sH_\rho(v)$ such that $\rho(D) < \rho$, then $v$ is a
minimal $L_\rho$--s.f.  \endproclaim

For example, $v\equiv 0$ is a minimal $L_\rho$--s.f.~because there
cannot be a negative $L_\rho$--s.f.~in $\Bbb T^2_P$ because of Proposition 8.6(2).

\demo{Proof} We first note that if $\rho(D) < \rho$, an
$L_\rho$--s.f. $v$ cannot be negative in all of $D$.  Indeed, let $q >
0$ solve (0.15) in a domain $D_1\subset\subset D$ such that
$\rho(D_1)=\rho$.  Choose $C=\min\limits_{z\in D_1}v(z)/(-q(z))$ so
that $v(z_0)+Cq(z_0)=0$ for some $z_0\in D_1$ and $v(z)+C q(z)\leq 0$
in $D_1$, but this contradicts the maximum principle Theorem 8.6(2).

Now suppose the Theorem is false, and let $v'$ be an $L_\rho$--s.m.~of
$v-\var$.  Then the function $v'-v$ is an $L_\rho$--s.f.~which is
negative in $D$.  Theorem 9.22 follows.\qed  \edm

It is possible to produce sufficient conditions for nonminimality.
For example,

\proclaim{Proposition 9.23} The function $v$ is nonminimal if
$v(z)\geq c$ or $L_\rho v-c > 0$ for some positive $c$ for all $z\in
\Bbb T^2_P$.  \endproclaim

This Proposition follows since $v\equiv c$ is an $L_\rho$--s.f.  But a
complete characterization of minimal $L_\rho$--subfunctions remains
open (see \cite {8, Problem 16.9}).

Vladimir Azarin,

\ Dept.of Math.\&Statistics
Bar-Ilan Univ.,Ramat-Gan 52900,Israel

\ \ e-mail: {\tt azarin\@macs.biu.ac.il}

David Drasin,

\ Math.Dept.,Purdue Univ.,West-Lafayette, IN 47907
USA,

\ \ e-mail: {\tt drasin\@math.purdue.edu}

Pietro Poggi-Corradini,

\ Dept.of Math.,Kansas State Univ., Manhattan, KS
66506, USA

\ \ e-mail: {\tt pietro\@math.ksu.edu}

\bye